%% file: VMmaxfinal.tex
\documentclass[graybox,envcountsect]{svmult}

\usepackage{mathptmx}       
\usepackage{helvet}         
\usepackage{courier}        
\usepackage{type1cm}        
\usepackage{makeidx}         
\usepackage{graphicx}        
\usepackage{multicol}        
\usepackage[bottom]{footmisc}

\usepackage{amsfonts,amsmath}

\makeindex             

\def \Z {{\mathbb Z}}
\def \R {{\mathbb R}}

\def \N {{\mathbb N}}

\def \PP {{\mathbb P}}
\def \EE {{\mathbb E}}

\def \ES {{\mathbb E}}
\def \PS {{\mathbb P}}

\def \c {{\rm c}}

\def \cP {{\mathcal P}}

\def\one{\rlap{\mbox{\small\rm 1}}\kern.15em 1}
\newcommand{\ind}{\boldsymbol{1}}

\begin{document}

\title*{Parabolic Anderson model with voter catalysts: dichotomy in the behavior 
of Lyapunov exponents}
\titlerunning{PAM with voter catalysts: dichotomy in the behavior of Lyapunov exponents} 
\author{G.\ Maillard, T.\ Mountford and S.\ Sch\"opfer}
\institute{G.\ Maillard \at CMI-LATP, Universit\'e de Provence,
39 rue F. Joliot-Curie, F-13453 Marseille Cedex 13, France, \email{maillard@cmi.univ-mrs.fr},
and EURANDOM, P.O.\ Box 513, 5600 MB Eindhoven, The Netherlands
\and T.\ Mountford \at  Institut de Math\'ematiques, \'Ecole Polytechnique F\'ed\'erale,
Station 8, 1015 Lausanne, Switzerland, \email{thomas.mountford@epfl.ch}
\and S.\ Sch\" opfer \at Institut de Math\'ematiques, \'Ecole Polytechnique F\'ed\'erale,
Station 8, 1015 Lausanne, Switzerland, \email{samuel.schoepfer@epfl.ch}}
%
%
\maketitle

\abstract*{}

\abstract{We consider the parabolic Anderson model $\partial u/\partial t =
\kappa\Delta u + \gamma\xi u$ with $u\colon\, \Z^d\times \R^+\to \R^+$, where
$\kappa\in\R^+$ is the diffusion constant, $\Delta$ is the discrete Laplacian,
$\gamma\in\R^+$ is the coupling constant, and
$\xi\colon\,\Z^d\times \R^+\to\{0,1\}$ is the voter model starting from Bernoulli
product measure $\nu_{\rho}$ with density $\rho\in (0,1)$. The solution of this
equation describes the evolution of a ``reactant'' $u$ under the influence of a
``catalyst'' $\xi$.\newline\indent
In G\"artner, den Hollander and Maillard \cite{garholmai10} the behavior of the
\emph{annealed} Lyapunov exponents, i.e., the exponential growth rates of
the successive moments of $u$ w.r.t.\ $\xi$, was investigated. It was shown that
these exponents exhibit an interesting dependence on the dimension and on
the diffusion constant.\newline\indent
In the present paper we address some questions left open in \cite{garholmai10}
by considering specifically when the Lyapunov exponents are the a priori maximal
value in terms of strong transience of the Markov process underlying the voter model.}


\section{Introduction}
\label{S1}


\subsection{Model}
\label{S1.1}

The parabolic Anderson model (PAM) is the partial differential equation
\begin{equation}
\label{pA}
\frac{\partial}{\partial t}u(x,t) = \kappa\Delta u(x,t) + \gamma\xi(x,t)u(x,t), \qquad x\in\Z^d,\,t\geq 0,
\end{equation}
with $u$ a $\R^+$-valued field, $\kappa\in\R^+$ a diffusion
constant, $\Delta$ the discrete Laplacian, acting on $u$ as
\begin{equation*}
\label{}
\Delta u(x,t) = \sum_{{y\in\Z^d} \atop {y\sim x}} [u(y,t)-u(x,t)]
\end{equation*}
($y\sim x$ meaning that $y$ is a nearest neighbor of $x$), $\gamma\in\R^+$ a coupling
constant and
\begin{equation*}
\label{}
\xi = (\xi_t)_{t \geq 0} \quad\text{with}\quad \xi_t = \{\xi_{t}(x):=\xi(x,t) \colon\,x\in\Z^d\}
\end{equation*}
the \emph{Voter Model} (VM) taking values in $\{0,1\}^{\Z^d\times\R^+}$. As initial condition, we choose
\begin{equation}
\label{inicond}
u(x,0) = 1, \qquad x\in \Z^d.
\end{equation}

One can interpret (\ref{pA}) in terms of population dynamics. Consider
a system of two types of particles, $A$ ``catalyst'' and $B$ ``reactant'', subject to:
\begin{itemize}
\item
$A$-particles evolve autonomously according to the voter dynamics;
\item
$B$-particles perform independent random walks at rate $2d\kappa$ and split
into two at a rate that is equal to $\gamma$ times the number of $A$-particles
present at the same location;
\item
the initial configuration of $B$-particles is one particle everywhere.
\end{itemize}
Then $u(x,t)$ can be interpreted as the average number of $B$-particles at site
$x$ at time $t$ conditioned on the evolution of the $A$-particles.


\subsection{Voter Model}
\label{S1.2}

The VM is the Markov process on $\{0, 1\}^{\Z^d}$ with
generator $L$ acting on cylindrical functions $f$ as
\begin{equation*}
L f(\eta)
=\sum_{x\in\Z^{d}}p(x,y)\sum_{{y\in\Z^d} \atop {y\sim x}}\Big(f(\eta^{x,y})- f(\eta)\Big),
\end{equation*}
where $p\colon\Z^d\times\Z^d\to[0,1]$ is the transition kernel of an irreducible random
walk and
$\eta^{x,y}$ is the configuration
\begin{equation*}
\begin{cases}
\eta^{x,y} (z)= \eta (z)& \forall z\not=x,\\
\eta^{x, y} (x)= \eta(y).
\end{cases}
\end{equation*}
In words, $\xi(x,t)=1$ and $\xi(x,t)=0$ mean the presence and the absence of a particle
at site $x$ at time $t$, respectively. Under the VM dynamics, the presence and absence
of particles are imposed according to the random walk transition kernel $p(\cdot,\cdot)$.

The VM was introduced independently by Clifford and Sudbury \cite{clisud73} and
by Holley and Liggett \cite{hollig75}, where the basic results concerning equilibria were
shown.
Let $(S_{t})_{t\geq 0}$ be the Markov semigroup associated with $L$, 
$p^{(s)}(x,y) = (1/2)[p(x,y)+p(y,x)]$, $x,y\in\Z^d$, be the \emph{symmetrized}
transition kernel associated with $p(\cdot,\cdot)$, and $\mu_\rho$ the \emph{equilibrium
measure} with density $\rho\in(0,1)$. When $p^{(s)}(\cdot\,,\cdot)$ is \emph{recurrent} all 
equilibria are \emph{trivial}, i.e., of the form $\mu_\rho = (1-\rho)\delta_0+\rho\delta_1$, 
while when $p^{(s)}(\cdot\,,\cdot)$ is \emph{transient} there are also \emph{non-trivial} 
equilibria, i.e., ergodic measures $\mu_\rho$, different from the previous one, which are 
the unique shift-invariant and ergodic equilibrium with density $\rho\in (0,1)$.

For both cases we 
have
\begin{equation}
\label{VMerg1}
\mu S_{t} \to \mu_\rho
\qquad \mbox{ weakly as } t\to\infty
\end{equation}
for any starting measure $\mu$ that is stationary and
ergodic with density $\rho$ (see Liggett \cite{lig85}, Corollary~V.1.13).
This is in particular the case for our choice $\mu:=\nu_{\rho}$, the Bernoulli
product measure with density $\rho\in(0,1)$.

\subsection{Lyapunov exponents}
\label{S1.3}

\begin{figure}[h]

\vspace{.5cm}
\begin{center}
\setlength{\unitlength}{0.5cm}
\begin{picture}(20,8)(0,0)
\put(0,0){\vector(1,0){8}}
\put(0,0){\vector(0,1){7.25}}
\put(-0.75,-0.75){$0$}
\put(-0.75,4.8){$\gamma$}
\put(8.5,-0.25){$\kappa$}
{\thicklines
\qbezier(0,5)(3,5)(7,5)
}
\put(5,6.5){$1\leq d\leq 4$}
\put(0,5){\circle*{.2}}
\put(-1,8){$\lambda_p(\kappa)$}
\put(13,0){\vector(1,0){8}}
\put(13,0){\vector(0,1){7.25}}
\put(12.25,-0.75){$0$}
\put(12.25,4.8){$\gamma$}
\put(11.75,0.8){$\rho\gamma$}
\put(21.5,-0.25){$\kappa$}
\put(12,8){$\lambda_p(\kappa)$}
{\thicklines
\qbezier(13,2.5)(14,1.4)(20,1.2)
\qbezier(13,3)(13.4,2.5)(14,2.3)
\qbezier(13,3.5)(13.4,3)(14,2.8)
}
\qbezier[60](13,1)(16.5,1)(20,1)
\qbezier[60](13,5)(16.5,5)(20,5)
\put(19,6.5){$d\geq 5$}
\put(13,2.5){\circle*{.2}}
\put(13,3){\circle*{.2}}
\put(13,3.5){\circle*{.2}}
\put(11.2,2.3){{\scriptsize $p=1$}}
\put(11.2,2.8){{\scriptsize $p=2$}}
\put(11.2,3.3){{\scriptsize $p=3$}}
\end{picture}
\vspace{.5cm}
\caption{\small Dichotomy of the behavior of $\kappa\mapsto\lambda_p(\kappa)$ when
$p(\cdot\,,\cdot)$ has zero mean and finite variance.}
\label{fig-lambda01}
\end{center}
\end{figure}
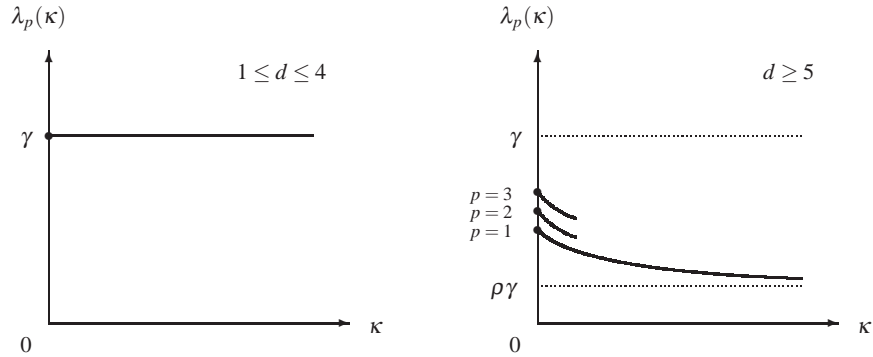

Our focus of interest will be on the $p$-th \emph{annealed} Lyapunov exponent,
defined by
\begin{equation}
\label{aLyapdef}
\lambda_p = \lim_{t\to\infty} \frac{1}{t} \log \EE_{\,\nu_{\rho}} \big([u(0,t)]^p \big)^{1/p},
\qquad p \in \N,
\end{equation}
which represents the exponential growth rate of the $p$-th moment of the
solution of the PAM (\ref{pA}), where $\EE_{\,\nu_{\rho}}$ denotes the expectation
w.r.t.\ the $\xi$-process starting from Bernoulli product measure $\nu_{\rho}$ with
density $\rho\in(0,1)$. Note that $\lambda_{p}$ depends on the parameters $\kappa$,
$d$, $\gamma$ and $\rho$ with the two latter being fixed from now. If the above limit
exists, then, by H\"older's inequality, $\kappa\mapsto\lambda_{p}(\kappa)$ satisfies
\begin{equation*}
\lambda_{p}(\kappa) \in [\rho\gamma,\gamma]
\qquad\forall \kappa\in[0,\infty).
\end{equation*}
The behavior of the annealed Lyapunov exponents with VM catalysts has already been
investigated by G\" artner, den Hollander and Maillard \cite{garholmai10}, where it was
shown that:
\begin{itemize}
\item
the Lyapunov exponents defined in (\ref{aLyapdef}) exist and do not depend on the
choice of the starting measure $\nu_{\rho} S_{T}$, $T\in[0,\infty]$, where
$\nu_{\rho}S_{\infty}:=\mu_{\rho}$ denotes the equilibrium measure of density $\rho$
(recall (\ref{VMerg1}));
\item
the  function $\kappa\mapsto \lambda_{p}(\kappa)$ is globally Lipschitz outside any
neighborhood of $0$ and satisfies $\lambda_{p}(\kappa)>\rho\gamma$ for all
$\kappa\in[0,\infty)$;
\item
the Lyapunov exponents satisfy
the following dichotomy (see Figure \ref{fig-lambda01}):
\begin{itemize}
\item
when $1\leq d\leq 4$,
if $p(\cdot,\cdot)$ has zero mean and finite variance, then
$\lambda_{p}(\kappa)=\gamma$ for all $\kappa\in[0,\infty)$;
\item
when $d\geq 5$,
\begin{itemize}
\item
$\lim_{\kappa\to 0}\lambda_{p}(\kappa)=\lambda_{p}(0)$;
\item
$\lim_{\kappa\to \infty}\lambda_{p}(\kappa)=\rho\gamma$;
\item
if $p(\cdot,\cdot)$ has zero mean and finite variance, then
$p\mapsto \lambda_{p}(\kappa)$ is strictly increasing for $\kappa\ll 1$.
\end{itemize}
\end{itemize}
\end{itemize}
The following questions were left open (see \cite{garholmai10}, Section 1.8):
\begin{description}[(Q1)]
\item[(Q1)]
Does $\lambda_{p}<\gamma$ when $d\geq 5$ if $p(\cdot,\cdot)$
has zero mean and finite variance?
\item[(Q2)]
Is there a full dichotomy in the behavior of the Lyapunov exponents? Namely,
$\lambda_{p}<\gamma$ if and only if $p^{(s)}(\cdot,\cdot)$ is strongly transient, i.e.,
\begin{equation*}
\label{strongtransient}
\int_{0}^\infty t p_{t}^{(s)}(0,0) \,dt<\infty \, .
\end{equation*}
\end{description}
Since any transition kernel $p(\cdot,\cdot)$ in $d\geq 5$ satisfies 
$\int_{0}^{\infty}t\, p_{t}(0,0)dt<\infty$, a positive answer to (Q2) will also ensure a positive 
one to (Q1) in the particular case when $p(\cdot,\cdot)$ is symmetric.
Theorems \ref{th1}--\ref{th3} in Section \ref{S1.4} give answers to question (Q2), depending
on the symmetry of $p(\cdot,\cdot)$. A positive answer to (Q1), given in Theorem \ref{th0},
can also be deduced from our proof of Theorem \ref{th1}.


By the Feynman-Kac formula, the solution of (\ref{pA}--\ref{inicond}) reads
\begin{equation*}
u(x,t) = \ES_{\,x}\left(\exp\left[\gamma\int_0^t \xi\left(X^\kappa(s),t-s\right)\,ds\right]\right),
\end{equation*}
where $X^\kappa=(X^\kappa(t))_{t \geq 0}$ is a simple random walk on $\Z^d$ with
step rate $2d\kappa$ and $\ES_{\,x}$ denotes the expectation with respect to $X^\kappa$
given $X^\kappa(0)=x$. This leads to the following representation of the Lyapunov exponents
\begin{equation*}
\lambda_{p}=\lim_{t\to \infty}\Lambda_{p}(t)
\end{equation*}
with
\begin{equation*}
\label{}
\Lambda_p(t)
= \frac{1}{pt} \log \big(\EE_{\,\nu_{\rho}}\otimes\ES_{\,0}^{\otimes p}\big)
\Bigg(\exp\Bigg[\gamma\int_0^t \sum_{j=1}^{p}
\xi\big(X_{j}^{\kappa}(s),t-s\big)\,ds\Bigg]\Bigg),
\end{equation*}
where $X_{j}^\kappa$, $j=1,\ldots,p$, are $p$ independent copies of $X^\kappa$.
In the above expression, the $\xi$ and $X^\kappa$ processes are evolving in time
reversed directions. It is nevertheless possible to let them run in the same
time evolution by using the following arguments. Let $\widetilde\Lambda_p(t)$
denote the $\xi$-\emph{time-reversal} analogue of $\Lambda_p(t)$ defined by
\begin{equation*}
\label{}
\widetilde\Lambda_p(t)=
\frac{1}{pt} \log \big(\EE_{\,\nu_{\rho}}\otimes\ES_{\,0}^{\otimes p}\big)
\Bigg(\exp\Bigg[\gamma\int_0^t \sum_{j=1}^{p}
\xi\big(X_{j}^{\kappa}(s),s\big)\,ds\Bigg]\Bigg)
\end{equation*}
and denote $\underline\Lambda_p(t) =$
\begin{eqnarray*}
&&\frac{1}{pt} \log \max_{x\in\Z^d}\big(\EE_{\,\nu_{\rho}}\otimes\ES_{\,0}^{\otimes p}\big)
\Bigg(\exp\Bigg[\gamma\int_0^t \sum_{j=1}^{p}
\xi\big(X_{j}^{\kappa}(s),t-s\big)\,ds\Bigg]
\prod_{j=1}^{p}\delta_{x}\big(X_{j}^\kappa(t)\big)\Bigg)\nonumber\\
&&=\frac{1}{pt} \log \max_{x\in\Z^d}\big(\EE_{\,\nu_{\rho}}\otimes\ES_{\,0}^{\otimes p}\big)
\Bigg(\exp\Bigg[\gamma\int_0^t \sum_{j=1}^{p}
\xi\big(X_{j}^{\kappa}(s),s\big)\,ds\Bigg]
\prod_{j=1}^{p}\delta_{x}\big(X_{j}^\kappa(t)\big)\Bigg),\\
\end{eqnarray*}
where in the last line we reverse the time of the $\xi$-process by using
that $\nu_{\rho}$ is shift-invariant and $X_{j}^\kappa$, $j=1,\ldots,p$, are
time-reversible.
As noted in \cite{garhol06}, Section 2.1,
$\lim_{t\to\infty} [\Lambda_{p}(t)-\underline\Lambda_{p}(t)]=0$
and, using the same argument,
$\lim_{t\to\infty} [\widetilde\Lambda_{p}(t)-\underline\Lambda_{p}(t)]=0$,
after which we can conclude that
\begin{equation*}
\label{fey-kac1}
\lambda_p(\kappa)=
\lim_{t\to\infty}\frac{1}{pt} \log \big(\EE_{\,\nu_{\rho}}\otimes\ES_{\,0}^{\otimes p}\big)
\Bigg(\exp\Bigg[\gamma\int_0^t \sum_{j=1}^{p}
\xi\big(X_{j}^{\kappa}(s),s\big)\,ds\Bigg]\Bigg).
\end{equation*}


\subsection{Main results}
\label{S1.4}

In what follows we give answers to questions (Q1) and (Q2) addressed
in \cite{garholmai10} concerning when the Lyapunov exponents are trivial,
i.e., equal to their a priori maximal value $\gamma$.

Our first theorem gives a positive answer to (Q1). It will be proved in
Section \ref{S2} as a consequence of the proof of Theorem \ref{th1}.

\begin{theorem}
\label{th0}
If $d\geq 5$ and $p(\cdot,\cdot)$ has zero mean and finite variance, then
$\lambda_{p}(\kappa)<\gamma$ for all $p\geq 1$ and $\kappa\in[0,\infty)$.
\end{theorem}

\begin{figure}[h]

\vspace{.5cm}
\begin{center}
\setlength{\unitlength}{0.5cm}
\begin{picture}(20,8)(0.5,0)
\put(0,0){\vector(1,0){8}}
\put(0,0){\vector(0,1){7.25}}
{\thicklines
\qbezier(0,6)(3,6)(7,6)
}
\put(-0.75,-0.75){$0$}
\put(-0.75,5.8){$\gamma$}
\put(8.5,-0.25){$\kappa$}
\put(-1,8){$\lambda_p(\kappa)$}
\put(4,7){{\scriptsize not strongly transient}}
\put(13,0){\vector(1,0){8}}
\put(13,0){\vector(0,1){7.25}}
{\thicklines
\qbezier(13,3.5)(14,1.4)(20,1.3)
}
\qbezier[60](13,1)(16.5,1)(20,1)
\qbezier[60](13,6)(16.5,6)(20,6)
\put(13,3.5){\circle*{.2}}
\put(12.25,-0.75){$0$}
\put(12.25,5.8){$\gamma$}
\put(11.75,0.8){$\rho\gamma$}
\put(21.5,-0.25){$\kappa$}
\put(12,8){$\lambda_p(\kappa)$}
\put(17.5,7){{\scriptsize strongly transient}}
\put(0,6){\circle*{.2}}
\end{picture}
\vspace{.5cm}
\caption{\small Full dichotomy of the behavior of $\kappa\mapsto\lambda_p(\kappa)$
when $p(\cdot,\cdot)$ is symmetric.}
\label{fig-lambda02}
\end{center}
\end{figure}
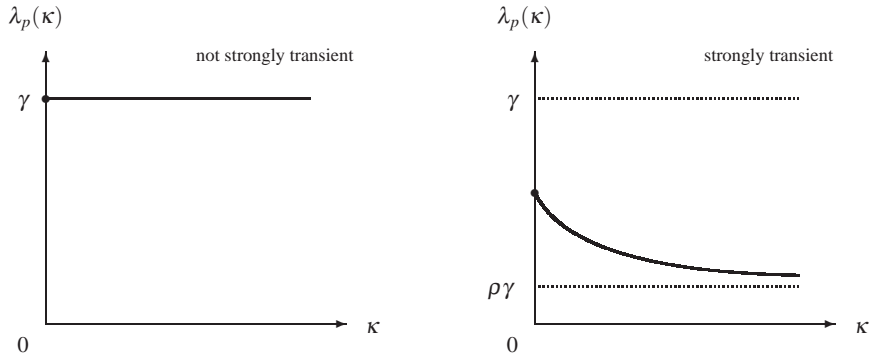

Our two next theorems state that the full dichotomy in (Q2) holds in the case
when $p(\cdot,\cdot)$ is symmetric (see Fig.\ \ref{fig-lambda02}). They will 
be proved in Section \ref{S2} and \ref{S3}, respectively.

\begin{theorem}
\label{th1}
If $p(\cdot,\cdot)$ is symmetric and strongly transient, then
$\lambda_{p}(\kappa)<\gamma$ for all $p\geq 1$ and $\kappa\in[0,\infty)$.
\end{theorem}

\begin{theorem}
\label{th2}
If $p(\cdot,\cdot)$ is symmetric and not strongly transient, then
$\lambda_{p}(\kappa)=\gamma$ for all $p\geq 1$ and $\kappa\in[0,\infty)$.
\end{theorem}


A similar full dichotomy also holds for the case where $\xi$ is symmetric exclusion
process in equilibrium, between recurrent and transient $p(\cdot,\cdot)$ (see
\cite{garholmai07}).

Our fourth theorem shows that this full dichotomy only holds for symmetric transition
kernels $p(\cdot,\cdot)$, ensuring that the assertion in (Q2) is not true in its full generality.

\begin{theorem}
\label{th3}
There exists $p(\cdot,\cdot)$ not symmetric with $p^{(s)}(\cdot,\cdot)$
not strongly transient such that
$\lambda_{p}(\kappa)<\gamma$ for all $p\geq 1$ and $\kappa\in[0,\infty)$.
\end{theorem}

 In the strongly transient regime, the following problems remain open:
\begin{description}[(a)]
\item[(a)] $\lim_{\kappa\to 0} \lambda_{p}(\kappa)=\lambda_{p}(0)$;
\item[(b)]
$\lim_{\kappa\to\infty}\lambda_{p}(\kappa)=\rho\gamma$;
\item[(c)]
$p\mapsto\lambda_{p}(\kappa)$ is strictly increasing for $\kappa\ll 1$;
\item[(d)]
$\kappa\mapsto \lambda_{p}(\kappa)$ is convex on $[0,\infty)$.
\end{description}
In \cite{garholmai10}, (a) and (b) were established when $d\geq 5$, and (c)
when $d\geq 5$ and $p(\cdot,\cdot)$ has zero mean and finite variance.
Their extension to the case when $p(\cdot,\cdot)$ is strongly transient
remains open.


In what follows, we use generic notation $\PP$ and $\EE$ for probability and 
expectation whatever the corresponding process is (even for joint processes) 
and denote $\xi_{s}(x):=\xi(s,x)$.


\section{Proof of Theorems \ref{th0} and \ref{th1}}
\label{S2}

We first give the proof of Theorem \ref{th1}. Recall that the transition kernel 
associated to the Voter Model $\xi$ is assumed to be symmetric. 
At the end of the section we will explain how to derive the proof 
of Theorem \ref{th0}.

We have to show that $\lambda_{p}(\kappa)<\gamma$ for all $\kappa\in[0,\infty)$.
In what follows we assume without loss of generality that $p=1$, the extension to
arbitrary $p\geq 1$ being straightforward.
Our approach is to pick a \emph{bad environment set} $B_{E}$ associated to the
$\xi$-process and a \emph{bad random walk set} $B_{W}$ associated to the random
walk $X^{\kappa}$ so that, for all $n \in \N$,
\begin{eqnarray}
\label{th1.1}
&&\EE\Big(\exp\Big[\gamma\int_{0}^n \xi_{s}(X^\kappa(s))\, ds\Big]\Big)\\
&&\qquad\leq \Big(\PP(B_{E})+\PS(B_{W})\Big)e^{\gamma n}
+ \EE
\Big(\one_{B_{E}^\c\cap B_{W}^\c}
\exp\Big[\gamma\int_{0}^n \xi_{s}(X^\kappa(s))\, ds\Big]\Big)
\nonumber
\end{eqnarray}
with, for some $0<\delta<1$,
\begin{equation}
\label{th1.2}
\PP(B_{E})\leq e^{-\delta n},
\qquad
\PS(B_{W})\leq e^{-\delta n},
\end{equation}
and,
\begin{equation}
\label{th1.3}
\int_{0}^{n} \xi_{s}(X^\kappa(s))\,ds \leq n(1-\delta)
\qquad\text{on } B_{E}^\c\cap B_{W}^\c.
\end{equation}
Since, combining (\ref{th1.1}--\ref{th1.3}), we obtain
\begin{equation*}
\label{th1.4}
\lim_{n\to\infty}\frac{1}{n}\log
\EE\Big(\exp\Big[\gamma\int_{0}^n \xi_{s}(X^\kappa(s))\, ds\Big]\Big)
< \gamma \, ,
\end{equation*}
it is enough to prove (\ref{th1.2}) and (\ref{th1.3}).

The proof of (\ref{th1.2}) is given in Sections \ref{S2.1}--\ref{S2.3} below,
and (\ref{th1.3}) will be obvious from our definitions of $B_E$ and $B_W$.

\subsection{Coarse-graining and skeletons}
\label{S2.1}

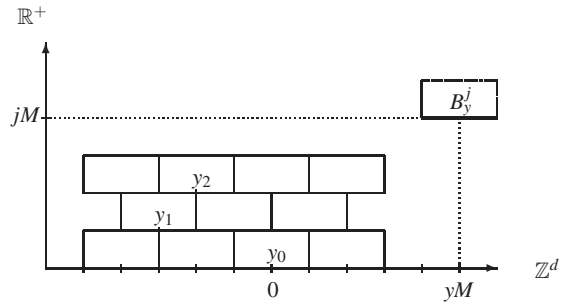
\begin{figure}[h]

\vspace{.5cm}
\begin{center}
\setlength{\unitlength}{0.5cm}
\begin{picture}(14,7)(0.5,0)
\put(0,0){\vector(1,0){12}}
\put(0,0){\vector(0,1){6}}
\put(1,-0.1){\line(0,1){1.1}}
\put(2,-0.1){\line(0,1){0.2}}
\put(3,-0.1){\line(0,1){1.1}}
\put(4,-0.1){\line(0,1){0.2}}
\put(5,-0.1){\line(0,1){1.1}}
\put(6,-0.1){\line(0,1){0.2}}
\put(7,-0.1){\line(0,1){1.1}}
\put(8,-0.1){\line(0,1){0.2}}
\put(9,-0.1){\line(0,1){1.1}}
\put(10,-0.1){\line(0,1){0.2}}
\put(11,-0.1){\line(0,1){0.2}}
\put(1,1){\line(1,0){8}}
\put(2,1){\line(0,1){1}}
\put(4,1){\line(0,1){1}}
\put(6,1){\line(0,1){1}}
\put(8,1){\line(0,1){1}}
\put(1,2){\line(1,0){8}}
\put(1,2){\line(0,1){1}}
\put(3,2){\line(0,1){1}}
\put(5,2){\line(0,1){1}}
\put(7,2){\line(0,1){1}}
\put(9,2){\line(0,1){1}}
\put(1,3){\line(1,0){8}}
{\thicklines
\put(10,4){\line(1,0){2}}
\put(10,4){\line(0,1){1}}
}
\qbezier[15](12,4)(12,4.5)(12,5)
\qbezier[30](10,5)(11,5)(12,5)
\put(10.75,4.25){$B_y^j$}
\qbezier[30](11,0)(11,2)(11,4)
\put(10.6,-0.75){$yM$}
\qbezier[80](0,4)(5,4)(10,4)
\put(-0.1,4){\line(1,0){0.2}}
\put(-0.95,3.9){$jM$}
\put(5.9,-0.75){$0$}
\put(5.9,0.25){$y_{0}$}
\put(2.9,1.25){$y_{1}$}
\put(3.9,2.25){$y_{2}$}
\put(3,0.9){\line(0,1){0.2}}
\put(4,1.9){\line(0,1){0.2}}
\put(-0.75,6.5){$\R^+$}
\put(13,-0.25){$\Z^d$}
\end{picture}
\vspace{.5cm}
\caption{\small Illustration of a $M$-skeleton $(y_{0},y_{1},y_{2},\ldots)\in \Xi$ 
and coarse-grained $B_y^j$ blocks.}
\label{fig-skeleton}
\end{center}
\end{figure}

Write $\Z_{\rm e}^d=2\Z^d$ and $\Z_{\rm o}^d=2\Z^d+1$, where $1=(1,\ldots,1)\in\Z^d$. 
We are going to use a \emph{coarse-graining representation} defined by a 
\emph{space-time block partition} $B_y^j$ and a \emph{random walk skeleton} 
$(y_i)_{i\geq 0}$. To that aim, for a fixed $M$, consider
\begin{equation*}
B_y^j = \prod_{k=1}^d \big[(y_{k}-1)M, (y_{k}+1)M\big) \times \big[jM,(j+1)M\big)
 \subset \Z^d \times \R_+\, ,
\end{equation*}
where $j \in \N_0:=\N \cup \{0\}$ and
\begin{equation*}
y \in \left\{\begin{array}{ll}\Z_{\rm e}^d & \text{when $j$ is even,} \\ 
\Z_{\rm o}^d & \text{when $j$ is odd.}
\end{array}\right.
\end{equation*}
Without loss of generality we can consider random walks trajectories on interval $[0,n]$
with $n \in \N$ multiple of $M$. Define the \emph{$M$-skeleton set} set by
\begin{equation*}
\Xi=\left\{\left(y_0,\ldots,y_{\frac{n}{M}}\right) \in (\Z^d)^{\frac{n}{M}+1} \colon
y_{2k} \in \Z_{\rm e}^d,\,\,\, y_{2k+1} \in \Z_{\rm o}^d\,\,\, \forall k \in \N_0\right\}
\end{equation*}
and the \emph{$M$-skeleton set associated to a random walk} $X$ by
\begin{equation*}
\Xi(X) = \left\{\left(y_0, \ldots ,y_{\frac{n}{M}}\right) \in \Xi\colon
X(kM) \in B_{y_k}^k \, \forall k \in \left\{0, \ldots ,n/M\right\}\right\} \, .
\end{equation*}
In what follows, we will consider the $M$-skeleton $\Xi(X^\kappa)$, but, as $X^\kappa$ 
starts from $0 \in \Z^d$, the first point of our $M$-skeleton will always be $y_{0}:=0\in\Z^d$
(see Fig.\ \ref{fig-skeleton}).


In the next lemma we prove that the number of $M$-skeletons not
oscillating too much is at most exponential in $n/M$. For that, define
\begin{equation}
\label{2.1.01}
\Xi_A=\left\{\big(y_0, \ldots ,y_{\frac{n}{M}}\big)\in \Xi \colon 
\sum_{j=1}^{n/M} (\|y_j-y_{j-1}\|_\infty-1) \leq \frac{n}{Md}\right\} \, ,
\end{equation}
where $\|\cdot\|_\infty$ is the standard $l_\infty$ norm, the set of all $M$-skeletons
that are \emph{appropriate}.
\begin{lemma}
\label{appropriate}
There exists some universal constant $K \in (1,\infty)$ such that, for any $n,M \in \N$,
\begin{equation*}
|\Xi_A| \leq K^{n/M} \, .
\end{equation*}
\end{lemma}
\begin{proof}
For any fixed $y_1 \in \Z^d$ and $N\in\N_{0}$, let
\begin{equation*}
I(N)=\big|\big\{y_2 \in \Z^d \colon \|y_1-y_2\|_\infty -1 = N \big\}\big|
\end{equation*}
be the number of elements of $\Z^d$ on the boundary of the cube of size $2N+3$
centered at $y_1$. For any $N \in \N$, we have $I(N) = (2N+3)^d-(2N+1)^d$ and 
$I(0)=3^d-1$, therefore, for any $N \in \N_0$, 
\begin{equation}
\label{IBd}
I(N) \leq 3^d(N+1)^d. 
\end{equation}
Define, for any $N,k \in \N$,
\begin{equation}
\label{2.1.02}
I(N,k)=\left|\left\{(y_j)_{0 \leq j \leq k} \in (\Z^d)^{k+1}\colon
\sum_{j=1}^k (\|y_j-y_{j-1}\|_\infty -1) = N\right\}\right| 
\end{equation}
the number of sequences in $(\Z^d)^{k+1}$ having size $N$. 
By (\ref{IBd}) and (\ref{2.1.02}), we have
\begin{eqnarray}
\label{2.1.03}
I(N,k) &=& \sum_{{(N_1, \ldots ,N_k)\colon}\atop{\sum_{i=1}^k N_i = N}}
\left( \prod_{i=1}^k I(N_i)\right) \nonumber\\
&\leq& 3^{dk} \sum_{{(N_1, \ldots ,N_k)\colon}\atop{\sum_{i=1}^k N_i = N}}
\left( \prod_{i=1}^k (N_i+1)\right)^d\nonumber\\
&\leq& 3^{dk} \binom{k+N+1}{N} 2^{dN} \, ,
\end{eqnarray}
where, in the last line,  we used that
\begin{equation*}
\max_{{(N_1, \ldots ,N_k)\colon}\atop{\sum_{i=1}^k N_i = N}} \prod_{i=1}^k (N_i+1) 
= 2^N.
\end{equation*}
Using (\ref{2.1.03}) and the fact that $\Xi \subset \Z^d$, we obtain
\begin{eqnarray*}
|\Xi_{A}| &\leq& \sum_{N=0}^{\left\lfloor \frac{n}{Md} \right\rfloor} I\left(N,\frac{n}{M}\right)
\leq 3^{\frac{n}{M}} \sum_{N=0}^{\left\lfloor \frac{n}{Md} \right\rfloor}
\binom{\frac{n}{M}+N+1}{N} 2^{dN}\\
&\leq & 3^{\frac{n}{M}} \sum_{N=0}^{\left\lfloor \frac{n}{Md} \right\rfloor}
\binom{\frac{2n}{M}}{N} 2^{dN}
\leq 3^{\frac{n}{M}} \sum_{N=0}^{\frac{2n}{M}}
\binom{\frac{2n}{M}}{N} 2^{dN}\\
& = & 3^{\frac{n}{M}} (2^d+1)^{\frac{2n}{M}}\, ,
\end{eqnarray*}
which ends the proof of the lemma.
\qed
\end{proof}

\subsection{The bad environment set $B_E$}
\label{S2.2}

This section is devoted to the proof of the leftmost part of (\ref{th1.2}) for suitable 
set $B_E$ defined below.

We say that an environment $\xi$ is \emph{good} w.r.t.\ an $M$-skeleton
$(y_0, \cdots ,y_{\frac{n}{M}})$ if we have
\begin{equation*}
\label{2.2.01}
\left| \left\{0 \leq j < \frac{n}{M}\colon
\exists (x,jM) \in B_{y_j}^j \text{ s.t. } \xi_s(x)=0 \,\,\,
\forall s \in [jM,jM+1]\right\}\right| \geq \frac{n}{4M} \, .
\end{equation*}
Since we want the environment to be good w.r.t.\ all appropriate $M$-skeleton,
we define the \emph{bad environment set} as
\begin{equation*}
B_E = \{\exists \text{ an $M$-skeleton} \in \Xi_A \text{ s.t. } \xi \text{ is not good w.r.t.\ it}\}.
\end{equation*}
In the next lemma we prove that for any fixed $M$-skeleton, the probability that $\xi$ is
not good w.r.t.\ it is at most exponentially small in $n/M$.  
\begin{lemma}
\label{be.lem}
Take $(y_i)_{0 \leq i \leq \frac{n}{M}} \in \Xi$ an $M$-skeleton. For $M$ big enough,
we have
\begin{equation*}
\PP\left(\xi \text{ is not good w.r.t. } (y_i)_{0 \leq i \leq \frac{n}{M}}\right)
\leq (4K)^{-n/M} \, ,
\end{equation*}
where $K$ is the universal constant defined in Lemma \ref{appropriate}.
\end{lemma}
Therefore, combining Lemmas \ref{appropriate} and \ref{be.lem}, we get
\begin{equation*}
\PP(B_E) \leq 4^{-n/M}
\end{equation*}
for $M$ big enough, from which we obtain the leftmost part of (\ref{th1.2}).

Before proving Lemma \ref{be.lem}, we first give an auxiliary lemma. In order to
study the evolution of the VM, we consider, as usual, the dual process, namely, 
a coalescing random system that evolves backward in time. To that aim, define 
$(X^{x,t}(s))_{0\leq s\leq t}$ to be the random walk starting from $0$ at time $t$, 
(i.e., $X^{x,t}(0)=0$). From the graphical representation of the VM, we can write 
$\xi_{t}(x)=\xi_{0}(X^{x,t}(t))$, $x\in\Z^d$, and therefore the the VM process can be 
expressed in terms of its initial configuration and a system of coalescing random walks.
Two random walks $X^{x,s}$ and $X^{x^\prime,s^\prime}$ with $s^\prime<s$
meet if there exists $u\leq s^\prime$ such that 
$X^{x^\prime,s^\prime}(s^\prime-u)=X^{x,s}(s-u)$. It is therefore the same to say
that $X^{x,s}$ and $X^{x^\prime,s^\prime}$ with $s^\prime<s$ meet (in some appropriate
time interval) if there exists $t\geq 0$ in this interval (by letting $t= s^\prime-u$) such that  
\begin{equation*}
X^{x^\prime,s^\prime}(t)=X^{x,s}(t+s-s^\prime)\, .
\end{equation*}
For convenience we will adopt this notation in the rest of the section.
\begin{lemma}
\label{2.2.lem1}
Take two independent random walks $X^{x,s}$ and $X^{x^\prime,s^\prime}$ with 
$s^\prime<s$. Then the probability they ever meet is bounded above by
\begin{equation*}
\int_{s-s^\prime}^\infty p_t(0,0)\, dt \, .
\end{equation*}
\end{lemma}

\begin{proof}
Consider the random variable
\begin{equation*}
W=\int_0^\infty 1 \left\{X^{x,s}(t)=X^{x^\prime,s^\prime}(t+(s-s^\prime))\right\} \, dt\, .
\end{equation*}
By symmetry, its expectation satisfies
\begin{eqnarray*}
\EE(W) &=& \int_0^\infty \PP\left(X^{0,0}(2t+s-s^\prime)=x-x^\prime\right) \, dt\\
 &\leq& \int_0^\infty \PP\left(X^{0,0}(2t+s-s^\prime)=0\right) \, dt\\
 &=& \frac{1}{2} \int_{s-s^\prime}^\infty p_t(0,0) \, dt \, .
\end{eqnarray*}
Moreover, we have
\begin{equation*}
\EE(W \,|\, W>0) = \int_0^\infty p_{2t}(0,0)\, dt \geq \frac{1}{2}
\end{equation*}
and then, since $\EE(W)=\EE(W \,|\, W>0)\, \PP(W>0)$, it follows that
\begin{equation*}
\PP(W>0) \leq 2\EE(W)=\int_{s-s^\prime}^\infty p_t(0,0) \, dt \, .
\end{equation*}
\qed
\end{proof}

We are now ready to prove the Lemma \ref{be.lem}.

\begin{proof}
Recall that $\xi_s(x)=\xi_0(X^{x,s}(s))$, where $\xi_0$ is distributed by a product
Bernoulli law with density $\rho\in(0,1)$. We first consider any $M$-skeleton 
$(y_0, \ldots ,y_{n/M}) \in \Xi$ (even not appropriate). For each $0 \leq j < n/M$, 
we choose $R$ sites $(x_1^j,jM), \ldots ,$ $(x_R^j,jM) \in B_{y_j}^j$ such that
\begin{equation}
\label{cond2}
\PP\left(\exists 0 \leq k,k^\prime \leq R,\,\,\, k\neq k^\prime\colon 
X^{x_k^j,jM}(s)=X^{x_{k^\prime}^j,jM}(s)
\text{ for some } s \in [0,jM]\right) \leq \epsilon
\end{equation}
for $\epsilon \ll 1$ to be specified later (see Fig.\ \ref{fig-x_k}). Remark that we first fix 
$\epsilon$ and $R$ and then we choose $M$ large enough so we can find these $R$ 
sites. As we are in the strongly transient regime, we know that these points exist. If two 
such random walks hit each other, then we freeze all the random walks issuing from 
the corresponding block $j$.

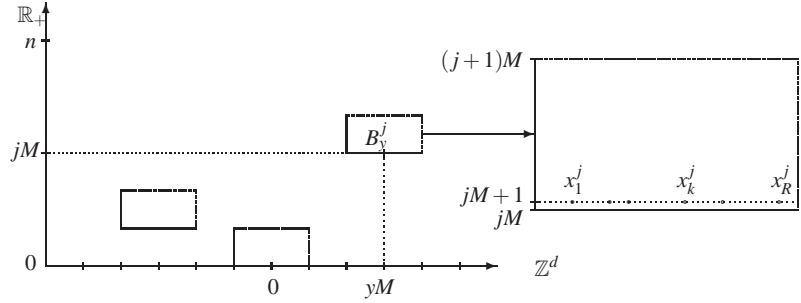
\begin{figure}[h]
\vspace{.5cm}
\begin{center}
\setlength{\unitlength}{0.5cm}
\begin{picture}(19,7)(0.5,0)
\put(0,0){\vector(1,0){12}}
\put(0,0){\vector(0,1){7}}
\put(1,-0.1){\line(0,1){0.2}}
\put(2,-0.1){\line(0,1){0.2}}
\put(3,-0.1){\line(0,1){0.2}}
\put(4,-0.1){\line(0,1){0.2}}
\put(5,-0.1){\line(0,1){0.2}}
\put(6,-0.1){\line(0,1){0.2}}
\put(7,-0.1){\line(0,1){0.2}}
\put(8,-0.1){\line(0,1){0.2}}
\put(9,-0.1){\line(0,1){0.2}}
\put(10,-0.1){\line(0,1){0.2}}
\put(11,-0.1){\line(0,1){0.2}}
\put(5,0){\line(0,1){1}}
\qbezier[15](7,0)(7,0.5)(7,1)
\qbezier[30](5,1)(6,1)(7,1)

\put(2,1){\line(1,0){2}}
\put(2,1){\line(0,1){1}}
\qbezier[15](4,1)(4,1.5)(4,2)
\qbezier[30](2,2)(3,2)(4,2)


\put(8,3){\line(1,0){2}}
\put(8,3){\line(0,1){1}}
\qbezier[15](10,3)(10,3.5)(10,4)
\qbezier[30](8,4)(9,4)(10,4)


\put(8.5,3.23){$B_y^j$}
\qbezier[20](9,0)(9,1.5)(9,3)
\put(9,2.9){\line(0,1){0.2}}
\put(8.5,-0.75){$yM$}
\qbezier[70](0,3)(4,3)(8,3)
\put(-0.1,3){\line(1,0){0.2}}
\put(-0.95,2.9){$jM$}
\put(-0.55,5.8){$n$}
\put(-0.1,6){\line(1,0){0.2}}
\put(-0.55,-0.1){$0$}
\put(5.9,-0.75){$0$}
\put(-0.75,6.5){$\R_+$}
\put(13,-0.25){$\Z^d$}

\put(10,3.5){\vector(1,0){3}}
\put(13,1.5){\line(1,0){7}}
\put(13,1.5){\line(0,1){4}}
\qbezier[60](20,1.5)(20,3.5)(20,5.5)
\qbezier[105](13,5.5)(16.5,5.5)(20,5.5)
\put(12.9,1.5){\line(1,0){0.2}}
\put(11.9,1.1){$jM$}
\put(12.9,1.7){\line(1,0){0.2}}
\put(11,1.7){$jM+1$}
\put(12.9,5.5){\line(1,0){0.2}}
\put(10.5,5.3){$(j+1)M$}
\qbezier[50](13,1.7)(16.5,1.7)(20,1.7)
\put(14,1.7){\circle{0.05}}
\put(13.8,2.1){$x_1^j$}
\put(15,1.7){\circle{0.05}}
\put(15.5,1.7){\circle{0.05}}
\put(17,1.7){\circle{0.05}}
\put(16.8,2.1){$x_k^j$}
\put(18,1.7){\circle{0.05}}
\put(19.5,1.7){\circle{0.05}}
\put(19.3,2.1){$x_R^j$}
\end{picture}
\vspace{.5cm}
\caption{\small Illustration of sites $(x_k^j,jM)$, $k\in\{1,\ldots,R\}$, 
for a fixed $M$-skeleton.}
\label{fig-x_k}
\end{center}
\end{figure}

For any $1 \leq j < \frac{n}{M}$, $1 \leq k \leq R$ for some $R>0$, we have
\begin{eqnarray*}
&&\EE\left(\sum_{j^\prime=j+1}^{\frac{n}{M}-1} \,\, \sum_{k^\prime=1}^R
1\left\{X^{x_{k^\prime}^{j^\prime},j^\prime M}(s+(j^\prime-j)M)
= X^{x_k^j,jM}(s) \text{ for some } s \in [0,jM]\right\}\right)\\
&&\qquad \leq R \sum_{j^\prime=j+1}^{\frac{n}{M}-1}
\int_{(j^\prime-j)M}^\infty p_t(0,0) \, dt\nonumber\\
&&\qquad \leq R \int_M^\infty\frac{t}{M}\,p_t(0,0) \, dt \, ,
\end{eqnarray*}
and therefore, summing over $1\leq k\leq R$, we get
\begin{eqnarray}
\label{2.2.08}
&&\EE\left(\sum_{j^\prime=j+1}^{\frac{n}{M}-1} \,\, \sum_{k,k^\prime=1}^R
1 \left\{X^{x_{k^\prime}^{j^\prime},j^\prime M}(s+(j^\prime-j)M)
= X^{x_k^j,jM}(s) \text{ for some } s \in [0,jM]\right\}\right)\nonumber\\
&&\qquad \leq \frac{R^2}{M} \int_M^\infty t \,p_t(0,0) \, dt \leq \epsilon^2,
\end{eqnarray}
for $M$ sufficiently large. Again, remark that we first fix $\epsilon$ and $R$, then 
we choose $M$ large enough.
For each $j$, we now define the filtration
\begin{equation*}
\label{2.2.10}
\mathcal{F}^j_t = \sigma\left(X^{x_k^j,jM}(s)\colon 0 \leq s \leq t,\,\, 1\leq k \leq R\right)
\end{equation*}
and the sub-martingale $Z^j(t):=$
\begin{equation*}
\EE\left(\sum_{j^\prime=j+1}^{\frac{n}{M}-1} \,\, \sum_{k,k^\prime=1}^R
1\left\{X^{x_{k^\prime}^{j^\prime},j^\prime M}(s+(j^\prime-j))
= X^{x_k^j,jM}(s) \text{ for some } s\in [0,t] \right\}\,\,\,\bigg|\,\,\, \mathcal{F}_t^j\right) \, ,
\end{equation*}
with the stopping time
\begin{equation*}
\tau^j=\inf\left\{t \geq 0\colon Z^j(t)>\epsilon \right\} \, .
\end{equation*}
We freeze every random walk issuing from block $j$ at time $jM \wedge \tau^j$.
Using (\ref{2.2.08}) and the Doob's inequality, we can see that
\begin{equation}
\label{2.2.12}
\PP(\tau^j \leq jM)\leq \PP\left(\sup_{0 \leq t \leq jM} Z^j_t \geq \epsilon\right) \leq \frac{R^2}{M\epsilon} \, \int_M^\infty t p_t(0,0) \, dt 
\leq \epsilon \, .
\end{equation}
Since, $Z^j$ is a continuous sub-martingale except at jump times of one of 
the random walks $X^{x_k^j,jM}$ and when a jump occurs, the increment is 
at most
\begin{equation*}
R \sum_{j^\prime=j+1}^{\frac{n}{M}-1} p_{(j^\prime-j)M}(0,0) \leq \epsilon^2
\end{equation*}
if $M$ is big enough. Therefore, for all $0 \leq t \leq \tau^j$, we get
\begin{equation}
\label{2.2.15}
Z^j(t) < \epsilon+\epsilon^2 \leq 2\epsilon \qquad 
\PP\text{-a.s.}\, .
\end{equation}
Now we say that $j$ is \emph{good} if
\begin{itemize}
\item $\tau^j> jM$;
\item the $R$ random walks $X^{x_1^j,jM}, \cdots ,X^{x_R^j,jM}$ do not meet;
\item the random walks $X^{x_k^j,jM}$ do not hit any point $x_{k^\prime}^{j^\prime}$
during interval $[(j-j^\prime )M-1,(j-j^\prime )M]$ for $j^\prime < j$;
\item the random walks $X^{x_k^j,jM}$ do not meet
$X^{x_{k^\prime}^{j^\prime},j^\prime M}$ for $j^\prime < j$.
\end{itemize}

By (\ref{2.2.12}), we know that the probability that the first condition does not occur
is smaller than $\epsilon$. By definition of the sites $x_1^j, \ldots,x_R^j$, we know
that the probability that random walks issuing from the same block $j$ at sites
$x_{k}^j$, $k=1,\ldots,R$, hit each other is smaller than $\epsilon$ (recall 
(\ref{cond2})). Moreover, the probability that the third condition does not occur is 
bounded from above by
\begin{equation*}
R \sum_{j^\prime=1}^{j-1} \, \int_{(j-j^\prime)M-1}^{(j-j^\prime)M} p_t(0,0) \, dt
\end{equation*}
which is as small as we want for $M$ large because we are in a transient case.
We still have to compute the probability that the fourth condition does not occur.
Furthermore, we can see that two random walks issuing from the same block evolve
independently until they meet, provided they do not meet a previous random walk.
From the above consideration, we get
\begin{eqnarray*}
\PP\left(j \text{ is not good}\,\,\big|\,\,\mathcal{G}^{j-1}\right)
&\leq& 3\epsilon + \sum_{j^\prime=1}^{j-1} \,\, \sum_{k,k^\prime=1}^R
\PP\left(X^{x_k^j,jM} \text{ meets } X^{x_{k^\prime}^{j^\prime},j^\prime M} \,\,\big|\,\,
\mathcal{G}^{j^\prime}\right) \, ,
\end{eqnarray*}
where 
\begin{equation*}
\label{2.2.16}
\mathcal{G}^j=
\sigma\left(X^{x_k^{j^\prime},j^\prime M}(s)\colon
1\leq k \leq R,\, 1\leq j^\prime \leq j,\, 0 \leq s \leq j^\prime M \wedge \tau^{j^\prime}\right)\, .
\end{equation*}
Here, we recall that $X^{x_k^j,jM}$ meets $X^{x_{k^\prime}^{j^\prime},j^\prime M}$
(with $j^\prime<j$) if we have
\begin{equation*}
X^{x_{k}^{j},j M}(s+(j-j^\prime)M)
= X^{x_k^{j^\prime},j^\prime M}(s) \text{ for some } 
s\in \big[0,\tau^{j^\prime} \wedge j^\prime M\big] \, .
\end{equation*}
By (\ref{2.2.15}), we have, for all $j^\prime$ fixed,
\begin{equation*}
\sum_{j=j^\prime+1}^{\frac{n}{M}-1} \,\, \sum_{k,k^\prime=1}^R
\PP\left(X^{x_k^j,jM} \text{ meets } X^{x_{k^\prime}^{j^\prime},j^\prime M}\,\,\,\bigg|\,\,\,
\mathcal{G}^{j^\prime}\right) \leq 2\epsilon \, .
\end{equation*}
Summing over all $1\leq j^\prime\leq n/M-2$, we get
\begin{equation*}
\sum_{j^\prime=1}^{\frac{n}{M}-2} \,\, \sum_{j=j^\prime+1}^{\frac{n}{M}-1} \,\,
\sum_{k,k^\prime=1}^R
\PP\left(X^{x_k^j,jM} \text{ meets } X^{x_{k^\prime}^{j^\prime},j^\prime M}\,\,\,\bigg|\,\,\,
\mathcal{G}^{j^\prime}\right) \leq \frac{n}{M} \, 2\epsilon \, ,
\end{equation*}
and then, interchanging the sums, we arrive at
\begin{equation*}
\sum_{j=2}^{\frac{n}{M}-1} \,\, \sum_{j^\prime=1}^{j-1} \,\,
\sum_{k,k^\prime=1}^R
\PP\left(X^{x_k^j,jM} \text{ meets } X^{x_{k^\prime}^{j^\prime},j^\prime M}\,\,\,\bigg|\,\,\,
\mathcal{G}^{j^\prime}\right) \leq \frac{n}{M} \, 2\epsilon \, .
\end{equation*}
Thus, there are at most $\left\lfloor \frac{n}{2M} \right\rfloor$ random positions $j$
with the property
\begin{equation*}
\sum_{j^\prime=1}^{j-1} \,\, \sum_{k,k^\prime=1}^R
\PP\left(X^{x_k^j,jM} \text{ meets } X^{x_{k^\prime}^{j^\prime},j^\prime M}\,\,\,\bigg|\,\,\,
\mathcal{G}^{j^\prime}\right) \geq 4\epsilon \, ,
\end{equation*}
and so at least $\left\lceil \frac{n}{2M} \right\rceil-2$ random positions $j$ have
the property
\begin{equation*}
\sum_{j^\prime=1}^{j-1} \,\, \sum_{k,k^\prime=1}^R
\PP\left(X^{x_k^j,jM} \text{ meets } X^{x_{k^\prime}^{j^\prime},j^\prime M}\,\,\,\bigg|\,\,\,
\mathcal{G}^{j^\prime}\right) < 4\epsilon \, .
\end{equation*}
For these random positions $j$, we then have
\begin{equation*}
\PP\big(j \text{ is good } |\,\, \mathcal{G}^{j-1}\big)
\geq 1-7\epsilon.
\end{equation*}
Using an elementary coupling, we have at least $\frac{n}{3M}$ positions that are
good with probability bounded by
\begin{equation*}
\PP\left(Y \geq \frac{n}{3M}\right) \geq 1-e^{-c(\epsilon)n/M}
\end{equation*}
for
\begin{equation*}
Y \sim B\left(\left\lceil \frac{n}{2M} \right\rceil-2,1-7\epsilon\right)
\end{equation*}
and $c(\epsilon) \to \infty$ as $\epsilon \to 0$.
Therefore, outside of a small probability $e^{-c(\epsilon)n/M}$, for at least $\frac{n}{3M}$ 
positions $j$, we have that the random walks $X^{x_k^{j},j M}$ are disjoint and so the 
values $\xi_{0}(X^{x_k^{j},j M}(s))$ are independent until time $s\leq j M$.
Then, using the fact that $(\xi_0(x))_{x \in \Z^d}$ are i.i.d.\ Bernoulli product  with
parameter $\rho$, we have that the number of positions $j$ so that there exists
$(x,j M) \in B^{j}_{y_{j}}$ with $\xi_s(x)=0$ and $s\in[jM, jM+1]$ is at least 
$\frac{n}{4M}$ outside the probability
\begin{equation*}
\PP\left(Y^\prime \geq \frac{n}{12M}\right) \leq \left(\frac{1}{4K}\right)^{n/M}
\end{equation*}
with $Y^\prime \sim B\left(\frac{n}{3M},\left((1-e^{-1}(1-\rho))^R\right)\right)$,
where $\epsilon$ and $R$ are chosen small and large enough, respectively.
\qed
\end{proof}

The proof of the leftmost part of (\ref{th1.2}) is now completed.

\subsection{The bad random walk set $B_W$}
\label{S2.3}

This section is devoted to the proof of the rightmost part of (\ref{th1.2}).

We are now interested in the random walk $X^\kappa$. We are going to prove that 
$(X^\kappa(s))_{0\leq s\leq n}$ has an appropriate $M$-skeleton and touches 
enough zeros outside a probability event exponentially small in $n$ (see Lemmas
\ref{bw1} and \ref{bw2} below, respectively). To define the bad random set $B_{W}$ 
announced in (\ref{}), we are going to define $B_{W}=B_{W_{1}}\cup B_{W_{2}}$, 
where the bad sets $B_{W_{1}}$ and $B_{W_{2}}$ correspond, respectively, to 
random walks trajectories $X^\kappa$ which do not have appropriate $M$-skeleton 
and do not touch enough sites occupied by a zero configuration of the VM. To be 
more precise, define 
\begin{equation*}
B_{W_{1}}=\big\{(\Xi(X^\kappa) \notin \Xi_A\big\}\, .
\end{equation*}   
In the next lemma, we prove that the probability of $B_{W_{1}}$ is exponentially small 
in $n$.
\begin{lemma}
\label{bw1}
Take $(X^\kappa(s))_{0 \leq s \leq n}$ and $\Xi(X^\kappa)=(y_0, \cdots ,y_{n/M})$ the
associated $M$-skeleton. Then, there exists a constant $K^\prime$ not depending on
$n$ such that
\begin{equation*}
\label{2.3.10}
\PP(B_{W_{1}})\leq e^{-K^\prime n} \, .
\end{equation*}
\end{lemma}

\begin{proof}

In order to have the random walk moving from one block of the skeleton to a 
nonadjacent one, the random walk has to make at least $M$ steps in the same 
direction. Keeping that in mind, define
\begin{equation*}
Y_j(s) = X^\kappa(jM+s)-X^\kappa(jM)
\end{equation*}
and let
\begin{equation*}
\tau_1^j
= \inf\{s\colon \|Y_j(s)\|_\infty \geq M\}, 
\qquad \tau_i^j
= \inf\big\{s > \tau_{i-1}^j\colon \|Y_j(s)-Y_j(\tau_{i-1}^j)\|_\infty \geq M\big\}.
\end{equation*}
Next, define
\begin{equation*}
W_i^j = \ind_{\{\tau^i_j < M\}}
\end{equation*}
and use an elementary coupling to have
\begin{equation*}
\PP(W_1^j=1) \leq e^{-c/M} \quad \text{and} \quad \PP(W_i^j=1|W_{i-1}^j=1) 
\leq e^{-c(i)/M} \leq e^{-c/M}
\end{equation*}
for some constants $c$ and $c(i)$ which verify $c(i) \geq c$. Therefore, we have 
that the number of jumps for the $j$th block is bounded above by the number of 
$W_i^j$ equals to 1. Using a coupling we can see that this is bounded above by 
a geometric law with parameter $e^{-c/M}$. Now if we consider all the blocks, by 
elementary properties of geometric random variables, we have
\begin{equation*}
\PP(B_{W_1}) \leq \PP\left(Y \geq \frac{n}{Md}\right) \leq e^{-c^\prime n}
\end{equation*}
for $Y \sim B\left(\frac{n}{M}\left(1+\frac{1}{d}\right),e^{-c/M}\right)$, 
some constant $c^\prime>0$, $M$ being large and the proof is done.
\qed
\end{proof}

Lemma \ref{bw1} proves the first part of the rightmost part of (\ref{th1.2}), namely the part concerned
with bad set $B_{W_{1}}$. Now we look at the number of times $X^\kappa$ stays on a
site where the VM has zero value. For that, define
\begin{equation*}
\tau_{i+1}=\inf\left\{t>\tau_i+1\colon \exists x \in \Z^d \text{ s.t. }
\|x-X^\kappa(t)\|_\infty \leq 2M,\, \xi_s(x)=0 \,\,\, \forall s \in [t,t+1] \right\}
\end{equation*}
with $\tau_0=0$ and
\begin{equation*}
k(M)=e^{-1/2}\inf_{\|x\|_\infty \leq 2M}\PP\big(X^\kappa(1/2)=x\,\,|\,\,X^\kappa(0)=0\big)
\end{equation*}
(remark that $k(M)$ does not depend on $n$ and is strictly positive).
Finally, we define
\begin{equation*}
\label{3.4.1}
B_{W_2}=\left\{(X^\kappa,\xi) \colon \tau_{\left\lfloor\frac{n}{2M}\right\rfloor}
\leq n-1 \text{ and } \int_0^n \xi_s(X^\kappa(s))\, ds 
\geq n \left(1-\frac{k(M)}{8M}\right)\right\} 
\end{equation*}
as being the bad set corresponding to random walks trajectories $X^\kappa$ 
which do not touch enough sites occupied by a zero configuration of the VM.
In the next lemma, we prove that such a set has an exponentially small probability
in $n$.
\begin{lemma}
\label{bw2}
There exists a constant $\delta > 0$ not depending on $n$, such that, for
$M$ big enough we have
\begin{equation*}
\PP(B_{W_2}) \leq e^{-n\delta}\, .
\end{equation*}
\end{lemma}

\begin{proof}
Take any realization of $\xi$ and for each time $\tau_i$, define a random variable
$Y_i$ which take value $1$ if $X^\kappa$ reaches a site with value zero at time
$\tau_i+\frac{1}{2}$ and stays at that point until time $\tau_i+1$, and takes value $0$
otherwise. Remark that after having fixed $n$, we can choose the state of $\xi_{t}$,
$t>n$, as we want, for example, full of zeros. Continue until 
$\tau_{\left\lfloor n/(2M) \right\rfloor}$ which is finite if $\xi$ is well chosen after time 
$n$. Using the strong Markov property, for every $k_i \in \{0,1\}$, we see that
\begin{equation*}
\PP(Y_i=1|Y_j=k_j, j<i) = \PP(Y_i=1|Y_{i-1}=k_{i-1}) \geq k(M) \, .
\end{equation*}
Then, it follows that $Y:=\sum_{i=1}^{\left\lfloor n/(2M) \right\rfloor} Y_i$ is stochastically 
greater than $Y^\prime$ the binomial random variable $B\left(\frac{n}{2M},k(M)\right)$. 
Moreover, if $\tau_{\left\lfloor n/(2M) \right\rfloor}\leq n-1$, we have
\begin{equation*}
\int_0^n \xi_s(X(s))\,ds \leq n-\frac{1}{2} \, Y \, .
\end{equation*}
Hence, we get
\begin{eqnarray*}
\PP(B_{W_2}) &\leq& \PP\left(\tau_{\left\lfloor \frac{n}{2M} \right\rfloor}
\leq n-1 \,\,\text{ and }\,\, n-\frac{1}{2}\,Y \geq n-\frac{nk(M)}{8M}\right)\\
&\leq& \PP\left(n-\frac{Y}{2} \geq n -\frac{nk(M)}{8M}\right)\nonumber\\
&=& \PP\left(Y \leq \frac{nk(M)}{4M}\right)\nonumber\\
&\leq& e^{-cn}
\end{eqnarray*}
for $n$ sufficiently large and $c$, a positive constant not depending on $n$.
This result being shown for any realization of $\xi$ (up to time $n$), this
ends the proof.
\qed
\end{proof}

Lemma \ref{bw2} proves the second part of the rightmost part of (\ref{th1.2}), namely 
the part concerned with bad set $B_{W_{2}}$. To complete the
proof of (\ref{th1.3}), it suffices to use the definition of $B_E$ and
$B_W=B_{W_1} \cup B_{W_2}$, and to remark that if $B_E$ and $B_{W_1}$ do not
occur, then the first condition of $B_{W_2}$, namely 
$\tau_{\lfloor n/2M \rfloor}\leq n-1$, is satisfied and therefore the second 
must be violated.

\subsection{Proof of Theorem \ref{th0}}
\label{S2.4}

The proof of Theorem \ref{th0} can be deduced from the proof of Theorem \ref{th1}.
Without assuming that $p(\cdot,\cdot)$ is symmetric, it is enough to see that
\begin{itemize}
\item
$\displaystyle\int_{0}^{\infty} t p_{t}(0,0)\, dt<\infty$, by local CLT, and
\item
there is enough symmetry because there exists some $C>0$ such that
$p_{t}(x,0)\leq C p_{t}(0,0)$ for all $x\in\Z^d$ and $t\in[0,\infty)$.
Therefore, Lemma \ref{be.lem} can still be applied.
\end{itemize}
From these two observations, the proof of Theorem \ref{th0} goes through
the same lines as the one of Theorem \ref{th1}.


\section{Proof of Theorem \ref{th2}}
\label{S3}

In this section we consider the Lyapunov exponents when the random walk kernel
associated to the voter model noise is symmetric and also not strongly transient, that is
\begin{equation*}
\int_{0}^{\infty} t p_{t} (0,0)\, dt = \infty \, .
\end{equation*}
We want to show that when $p(\cdot,\cdot)$ is
symmetric and not strongly transient, then
\begin{equation*}
\lambda_{p} (\kappa) \equiv \gamma
\qquad \forall \kappa \in (0, \infty),\,\, \forall p \geq 1 \, .
\end{equation*}
Since the result is easily seen for recurrent random walks, we can and will assume 
in the following that
\begin{equation*}
\int_{0}^{\infty} p_{t} (0,0)\, dt < \infty \, .
\end{equation*}

Given the reasoning of \cite{garholmai07}, Section 3.1 and \cite{garholmai10}, 
Section 5.1 this result will follow from Proposition $\ref{prop2}$, below. Consider, 
in the graphical representation associated to the VM $\xi$,
\begin{equation*}
\chi(t) : = \text{ number of distinct coalescing random walks produced on } \{0\} 
\times [0,t]
\end{equation*}
(This quantity is discussed in Bramson, Cox and Griffeath \cite{bramcoxgri88}).
\begin{proposition}
\label{prop2}
Assume that $p(\cdot,\cdot)$, is symmetric and not strongly transient,
then for any $\epsilon > 0$, we have that
\begin{equation*}
\lim_{t \to \infty} \PP \big(\chi(t) \leq \epsilon t\big) = 1 \, .
\end{equation*}
\end{proposition}

Before proving Proposition \ref{prop2}, we will first give the proof of Theorem \ref{th2}.
\begin{proof}
From the graphical representation of the VM and Proposition \ref{prop2}, we can see that 
for all $\delta > 0$ and $M < \infty$,
\begin{equation*}
\PP\Big(\xi (x, s) = 1 \,\,\, \forall \|x \|_{\infty} \leq M, \,\, \forall s \in [0, t]\Big)
\geq e^{-\delta t}
\end{equation*}
for all $t$ sufficiently large (see \cite{garholmai10}, proof of Lemma 5.1). Thus, 
just as in $\cite{garholmai10}$, Section 5.1, we have for all $p \geq 1$,
\begin{eqnarray*}
&&\EE ([u(0, t)]^{p})\\
&&\qquad\geq e^{\gamma p t} \PP\Big(\|X^\kappa(s) \|_{\infty} < M \,\,\, 
\forall s\in [0,t]\Big)^{p}\,
\PP\Big(\xi (x, s) = 1 \,\,\, \forall \| x \|_{\infty} < M,\, \forall s\in [0,t]\Big)\\
&&\qquad\geq e^{t (\gamma p - \delta -c(M) p)}
\end{eqnarray*}
for $c(M)\to 0$ as $M\to \infty$. From this, it is immediate that
\begin{equation*}
\lim_{t\to\infty}\frac{1}{t} \log \EE([u(0, t)]^{p})^{1/p} = \gamma.
\end{equation*}
\qed
\end{proof}

To prove Proposition \ref{prop2}, we consider the following system of coalescing
random walks
\begin{equation*}
I= \{X^{t} \colon t \in \cP \}
\end{equation*}
for $\cP$ a two sided, rate one Poisson process and $X^{t}$ a random walk defined
on $s \in [t, \infty)$, starting at $0$ at time $t$ (we could equally well consider a system 
of random walks indexed by $h \Z$ for some constant $h$). The coalescence is such 
that for $t < t^\prime \in \cP$, $X^{t}$, $X^{t^\prime}$ evolve independently until 
$T^{t, t^\prime}= \inf\{s > t^\prime \colon X^{t}(s) = X^{t^\prime}(s)\}$, and then, for 
$s \geq  T^{t, t^\prime}$, $X^{t}(s) = X^{t^\prime}(s)$.

We will be interested in the density or number of distinct random walks at certain times.
To aid this line we will adopt a labelling procedure for the random walks, whereby
effectively when two random walks meet for the first time, one of them (chosen at
random) dies; in this optic the number of distinct random walks will be the number
still alive.  Our labeling scheme involves defining for each $t\in \cP, $ the label
process $l_{s}^{t}$ for $s \geq t-$ (it will be helpful to be able to define  $l_{t-}^{t} = t$,
though since at time $t$ there may well be other random walks present at the origin,
it will not necessarily be the case that $ l_{t}^{t} = t$).  These processes will be defined
by the following properties:
\begin{itemize}
\item
if for  $t \ne t^\prime \in \cP, X^t(s) \ne X^{t^\prime}(s)$, then  
$l^{t }_s \ne l^{t^\prime}_{s}$;
\item
if $t_{1}, t_{2} , \ldots,t_{r} $ are elements of $\cP$, then at
$s \geq \max \{ t_{1}, \ldots, t_{r}   \}$, if
$X^{t_{1}}(s) =X^{t_{2}}(s) = \cdots= X^{t_{r}}(s)$,
then $l^{t_{1 } }_{s}=l^{t_{2}}_{s} = \cdots =l^{t_{r}}_{s} =u $
for some $ u \in \cP$ with $ X^{t_{1}}(s) =X^{u}(s)$;
\item
if for $t \ne t^\prime \in \cP$, $X^t $ meets $X^{t^\prime}$ for the first time at $s$,
then independently of past and future random walks or labeling
decisions $l^{t }_s = l^{t^\prime}_s = l^{t^\prime}_{s-}$ with probability $\frac12 $
and with equal probability $l^{t }_s = l^{t^\prime}_s = l^{t}_{s-}$;
\item
the process $l^t_{s}$ can only change at moments
where $X^t $ meets a distinct random walk for the first time.
\end{itemize}



For $t \in \cP$, $s > t$, we say that $t$ is \emph{alive at time} $s$, if $l^t_s = t$; it
\emph{dies} at time $s$ if $l_{s-}^{t} = t, l_{s}^{t}  \ne t$.   We say $X^t$, $X^u$
\emph{coalesce} at time $s$ if this is the first time at which the two labels are equal.
The following are easily seen:
\begin{itemize}
\item
the events $A_{s}^{t} = \{l_{s}^{t} = t \}$ for $t \in \cP$ are decreasing in $s$;
\item
$A_{s}^{t}$ depends only on the random motions of the coalescing random
walks and on the labeling choices involving $X^{t}$;
\item
for $s > 0,  $ the number of independent random walks $X^{t}(s) $,
$t \in \cP \cap [-n, 0]$ is simply equal to the number of distinct labels
$l_{s}^{t}$, $t \in \cP \cap [-n, 0]$.
\end{itemize}
Let
\begin{equation}
\label{Atslim}
c_0 = \lim_{s \to \infty} \PP^{t} (A_{s}^{t}) \in [0,1],
\end{equation}
according to palm measure, $\PP^{t}$, for $t \in \cP$. We obtain easily:
\begin{proposition}
\label{lem3}
\begin{equation*}
\lim_{s\to\infty}\frac{1}{s}
\big| \big\{ \text{distinct random walks } X^{t}(0) \colon t \in \cP \cap [-s,0) \big\}\big|
= c_0\qquad\text{a.s.}
\end{equation*}
\end{proposition}

\begin{proof}
Using the definition of $c_0$ in (\ref{Atslim}) and ergodicity of the system we see that the limit is greater than $c_0$. Then, Lemma \ref{prop4} gives the result.
\end{proof}

\begin{lemma}
\label{prop4}
For $c_0$ as defined in (\ref{Atslim}), for each $\epsilon > 0$, there exists 
$R < \infty$ so that if we consider the finite system of coalescing random walks
${(X^{t})}_{t \in (-R, 0] \cap \cP}$, then with probability at least $ 1-\epsilon$ at 
time $R$ there are less than $(c_0 + \epsilon)R$ distinct random walks labels.
\end{lemma}

\begin{proof}
By definition of $c_0$, for all $\epsilon>0$ there exists a $T_{0}$ so that
\begin{equation*}
\PP^{0} \big(\text{ label } 0 \text{ is alive at time } s\big) < c_0 + \frac{\epsilon}{100}
\qquad \forall s \geq T_{0}.
\end{equation*}
Now pick $R_{1}$ so that
\begin{equation*}
\PP \big( \|X^{0}(s)\|_{\infty} \leq R_{1} \,\,\, \forall s\in(0,T_{0})\big)
\geq 1- \frac{\epsilon}{100}.
\end{equation*}
Therefore,
\begin{eqnarray*}
&&\PP^{0} \big(\text{ label } 0 \text{ is not alive at time } s\geq T_{0}
\text{ and } \|X^{0}(s)\|_{\infty} \leq R_{1} \,\,\, \forall s\in(0,T_{0})\big) \\
&&\qquad\geq 1 - c_{0} - \frac{2\epsilon}{100}\, .
\end{eqnarray*}
We then pick $T_{1}$ so that
\begin{equation*}
\PP\big(\exists t \in \cP \cap {[-T_{1}, T_{1}]}^{\c} \colon 
\big\|X^{t}(s)\big\|_{\infty} \leq R_{1}\,\,\,
\text{for some } s \in (0, T_{0})\big) < \frac{\epsilon}{100}.
\end{equation*}
Thus
\begin{equation*}
\PP^0\big(\exists t \in \cP \cap [-T_{1}, T_{1}] \setminus \{0\} \colon 
l_{R_1}^{0}=t\big) \geq 1-c_0- \frac{\epsilon}{30}.
\end{equation*}
From the translation invariant property of the system and ergodicity if
\begin{equation*}
\lambda_{s} := \Big|\Big\{t \in [-s,s] \cap \cP\colon X^{t} \text{ loses its label 
to a random walk } X^{t^\prime} \text{ with } |t-t^\prime| \leq T_{1}\Big\}\Big|,
\end{equation*}
then
\begin{equation*}
\liminf_{s \rightarrow \infty} \frac{\lambda_{s}}{2 s} 
\geq 1-c_0-\frac{\epsilon}{30}\qquad \text{a.s.}
\end{equation*}
The result now follows easily.
\qed
\end{proof}
Proposition \ref{prop2} will be proven by showing:
\begin{proposition}
\label{prop3}
If $p(\cdot,\cdot)$ is symmetric and not strongly transient, then $c_0 =0$.
\end{proposition}

The proof of Proposition \ref{prop3} will work for any Poisson process rate, in
particular for $\cP$ having rate $M \gg 1$.  The distinct random walks treated in 
Proposition \ref{prop2} can be divided into those coalesced with a random walk from 
the system derived from $\cP$ (and so by Proposition \ref{prop3} of small ``density'') 
and those uncoalesced (also of small ``density'' if $M$ is large). Thus
Proposition \ref{prop2} follows almost immediately from Proposition \ref{prop3}.

The argument for Proposition \ref{prop3} is low level and intuitive. We argue by
contradiction and suppose
that $c_0 > 0$. From this we can deduce, loosely speaking, that after a certain
time either a random walk has lost its original label, or it will keep it forever. We
then introduce coupling on these random walks so that we may regard these
random walks as essentially independent random walks starting at $0$ (at
different times). We then introduce convenient comparison systems so that we
can analyze subsequent coalescences.  We will use automatically, without
reference, the following ``obvious" result:
\begin{lemma}
Consider two collections of coalescing random walks $\{Y^i\} $ and $\{(Y^\prime)^i\}$
for $i$ in some index set. If the coalescence rule is weaker for the  $\{(Y^\prime)^i\}$
system, in that if two walks $(Y^\prime)^i$ and $(Y^\prime)^j$ are permitted to
coalesce at time $t$, then so are $Y^{i }$ and $Y^{j}$, then there is a coupling of the
two systems so that the weaker contains the stronger.
\end{lemma}

We now fix $\epsilon > 0$ so that $\epsilon \ll c_0$ (by hypothesis $c_0 >0)$.
We choose  $R$ according to Lemma $\ref{prop4}$ and divide up time into
intervals $I_{j} = [jR, (j+1)R)$. We first consider the coalescing system where random
walks $X^{t}$, $X^{t^\prime}$, $t, t^\prime \in \cP$, can only ``coalesce'' (or destroy 
a label $t$ or $t^\prime$) if $t, t^\prime$ are in the same $I_{j}$ interval. Thus we 
have a system of random walks that is invariant to time shifts by integer multiples of 
$R$. We now introduce a system of random walks $Y^{t}$,
$t \in V := \cup_j \{ [jR, (j+1)R) \cap jR+\frac{1}{c_0} \Z \}$.
The random walks $Y^{t}$, $t \in \big[jR,(j+1)R\big)$ are not permitted to coalesce 
up until time $(j+1)R$ (at least) and will evolve independently of the system 
${(X^{t})}_{t \in \cP}$ until time $(j+1)R$. We will match up the points in $V \cap I_j$, with those in 
$\cP \cap I_{j}\cap K$ in a maximal measurable way for 
$K = \{t \in \cP \colon \text{label } t \text{ survives to time } (j+1)R\}$.

\begin{lemma}
\label{lem8}
Unmatched points in $\cup_{j} \cP \cap I_{j} \cap K$ and in $V$ have density
less than $2 \epsilon$ for $R$ fixed sufficiently large.
\end{lemma}

Remark: the system is not translation invariant with respect to all shifts but it 
possesses enough invariance for us to speak of densities.

We similarly have
\begin{lemma}
Unmatched $Y$ particles have density less than $2\epsilon$ for $R$ fixed 
sufficiently large.
\end{lemma}

It is elementary that two random walks $X$, $Z$ can be coupled so that for $t$ 
sufficiently large  $X(t) = Z(t)$. For given $\epsilon > 0$ we choose $M_{0}$ and 
then $M_{1}$ so that
\begin{equation}
\label{M0def}
\PP\Big(\sup_{t \leq R} \|X(t)\|_{\infty} \geq M_{0} \Big) < \frac{\epsilon}{10}
\end{equation}
and
\begin{equation}
\label{M1def}
\sup_{|z|\leq 2 M_{0}} \PP \Big(X_0, Z_{z} \text{ not coupled by time }
 M_{1}\Big) < \frac{\epsilon}{10},
\end{equation}
where $X_0$ and $Z_{z}$ denote that the random walks $X$ and $Z$ start at point 
$0$ and point $z$ respectively.

We then (on interval $I_{j}$) couple systems $Y$ and $X$ by letting married
pairs $Y^{t}$, $X^{t^\prime}$, $t \in V$, $t^\prime \in \cP \cap K$, evolve independently of 
other $Y$, $X$ random walks so that they couple by time $M_{1} + (j+1) R$ with 
probability at least $1- \frac{3\epsilon}{10}$.

Thus we have two types of random walk labels, $l^t$, for the $X$ system which are
equal to $t$ at time $t+R+M_{1}$: those for which the associated random walk was
paired with a $Y$ random walk and such that the random walks have coupled by time
$t+R+M_{1}$ said to be \emph{coupled} and the others, said to be \emph{decoupled}.
Similarly for the points in $V$ associated to $Y$ random walks.  We note that the 
foregoing implies that the  density of uncoupled labels is bounded by 
$2\epsilon+\frac{3\epsilon}{10} \leq 3 \epsilon$.  
The point is that modulo this small density, we have an identification of the coalescing 
$X$ random walks and the $Y$ random walks.

We now try to show that enough $Y$ particles will coalesce in a subsequent time interval 
to imply that there will be a significant decrease in surviving labels for the $X$ system. To 
do this we must bear in mind that, essentially, it will be sufficient to show a decrease in the 
density of $Y$ random walk labels definitely greater than $\epsilon$.  Secondly, as already 
noted, we will adopt a coalescence scheme that is a little complicated namely 
$Y^{i}, Y^{i^\prime}$ in $V$ can only ``coalesce'' at time $t \geq \max\{i, i^\prime\}$ if
\begin{itemize}
\item
$(t-i^\prime, t-i)$ are in some time set to be specified;
\item
for $i<i^\prime$, $\frac{t-i^\prime}{i^\prime-i} \in \left(\frac{9}{10}, \frac{11}{10} \right)$.
\end{itemize}
We now begin to specify our coalescence rules for the random walk system $\{Y^{i}\}_{i \in V}$.
The objective here will be to facilitate the necessary calculations. A first objective is to have
coalescence of $Y^{i}, Y^{i^\prime}$ at times $t > \max\{i, i^\prime\}$ so that $p_{ t}(0,0)$ is well behaved
around $t-i$, $t-i^\prime$. It follows from symmetry of the random walk that $t \mapsto p_{t}(0,0)$ is 
decreasing. The problem
we address is that it is not immediate how to achieve bounds in the opposite direction. This is
the purpose of the next result.
\begin{lemma}
\label{lem5}
Consider positive $\{a_{n}\}_{n\geq 0}$ so that $\sum_{n=0}^{\infty} a_{n} = \infty$.
For all $r \in \Z_{+}$, there exists a subsequence $\{a_{n_{i}}\}_{i\geq 0}$ so that\\
(i) $\sum a_{n_{i}} = \infty$;\\
(ii) $a_{n_{i}} > \frac{1}{2} a_{n_{i}-r}$, $\forall i\geq 0$.
\end{lemma}
\begin{proof}
If $r > 1$ we may consider the $r$ subsequences ${\{a_{ri +j}\}}_{i\geq 0}$ for
$j \in \{0, 1, \cdots, r-1\}$. At least one of these must satisfy $\sum {a_{ri+j}} = \infty$ so,
without loss of generality, we take $r=1$.

Now we classify $i$ as good or bad according to whether $a_{i} > a_{i-1}/2$ or not.
This decomposes $\Z$ into intervals of bad sites, alternating with intervals of good sites.
By geometric bounds, the sum of bad sites is bounded by the sum of the good $a_{i}$ for
which $i$ is the right end point of a good interval. Thus we have
\begin{equation*}
\sum_{i \textrm{ good}} a_{i}= \infty,
\end{equation*}
from which the result is immediate.
\qed
\end{proof}

\begin{corollary}
\label{lem5.1}
For our symmetric kernel $p_{t}(0,0)$ we can find $n_{i} \uparrow \infty$ so that\\
(i) $\displaystyle \sum_{i} \int_{2^{n_{i}}}^{2^{n_{i}+1}} t p_{t}(0,0) dt = \infty$;\\
(ii) $\displaystyle p_{2^{n_{i}-1}}(0,0) \leq 2^{12}p_{2 ^{n_{i} +3}} (0,0)$.
\end{corollary}
\begin{proof}
In Lemma \ref{lem5} take
\begin{equation*}
a_{n}= \int_{2^{n+3}}^{2^{n+4}} t p_{t} (0,0) dt
\end{equation*}
and take $r=5$.  Then by the monotonicity of $t \rightarrow p_t(0,0)$, we have
\begin{equation*}
2^{2n_{i}+7}p_{2^{n_{i}+3}} (0,0)  \geq a_{n_{i}}  \geq \frac{1}{2} \, a_{n_{i}-5}
\geq \frac{1}{2} \, 2^{2n_{i}-4}p_{2^{n_{i}-1}} (0,0).
\end{equation*}
\qed
\end{proof}
We fix such a sequence $\{n_{j}\}_{j\geq 1}$ once and for all.

We assume, as we may, that
$n_{j} < n_{j+1} -4$ for all $j\geq 1$ and also assume again, as we may, that
\begin{equation*}
\int_{2^{n_{j}}}^{2^{n_{j} +1}} t p_{t} (0,0) dt < \frac{\epsilon}{100}
\end{equation*}
for all $j\geq 1$.
We are now ready to consider our coalescence rules. We choose
$\epsilon \ll \alpha \ll 1$ (we will fully specify $\alpha$ later on but we feel it
more natural to defer the technical relations). We then choose $k_{0}$ so
that $2^{n_{k_{0}}} > R+M_{1}$ with $R$ as in Lemma $\ref{prop4}$ and 
Lemma $\ref{lem8}$ and $M_{1}$ as in (\ref{M1def}), and
\begin{equation}
\label{k1def}
k_{1} := \inf \bigg\{k > k_{0} \colon
\sum_{j=k_{0}}^{k} \int_{2^{n_{j}}}^{2^{n_{j}+1}} t p_{t} (0,0)\, dt > \alpha\bigg\}.
\end{equation}
We have coalescence between $Y^{i}$ and $Y^{i^\prime}$, for $i < i^\prime$ only at 
$t \in [i^\prime+2^{n_{j}},i^\prime+ 2^{n_{j}+1} ]$, $j \in [k_{0}, k_{1}]$ if
\begin{description}[(a)]
\item[(a)]
$\frac{t-i^\prime}{i^\prime-i} \in  (9/10, 11/10)$;
\item[(b)]
the interval of $t \in [i^\prime+2^{n_{j}}, i^\prime+2^{n_{j}+1} ]$ satisfying (a) is  of 
length at least $2$.
\end{description}
We say $(i,i^\prime)$ and $(i^\prime,i)$ are in $j$ and write $(i,i^\prime) \in j$ if the above 
relations hold. To show that sufficient coalescence occurs, we essentially use Bonferroni 
inequalities (see, e.g., \cite{durrett05}, p.\ 21).

To aid our argument we introduce a family of independent (non coalescing) random
walks $\{(Z^i(s))_{s \geq i}\}_{i \in V}$ such that for each $i \in V$, $Y^i(s) = Z^i(s)$ for
$s \geq i$ such that $l_s^i = i$.
In the following we will deal with random walks $Y^0, Z^0$, but lack of total translation
invariance notwithstanding,
it will be easy to see that all bounds obtained for these random walks remain valid for
more general random walks $Y^i, Z^i$.
For a given random walk $Y^{0} $, say, the probability that $Y^{0}$ is killed by
$Y^{i}$ (with $i$ possible in the sense of the above rules) is in principle a
complicated event given the whole system of coalescing random walks.
Certainly  the event
\begin{equation*}
\big\{Z^{0} \text{ meets }Z^{i} \text{ in appropriate time interval
after first having met }Z^{k}\big\}
\end{equation*}
is easier to deal with than the corresponding $Y$ event.
From this point on we will shorten our phraseology by taking ``$Z^i $ hits $Z^k$" to
mean that $Z^i $ meets  $Z^k$ at a time $t$ satisfying the conditions (a) and (b) above with
respect to $i$, $k$.

For ${(Z^{0}(s))}_{s \geq 0}$ and ${(Z^{i}(s))}_{s \geq i}$ independent random
walks each beginning at 0, we first estimate
\begin{equation*}
\sum_{i} \PP\big(Z^{0} \text{ hits } Z^{i} \big).
\end{equation*}
This of course decomposes as
\begin{equation*}
\sum_{j} \sum_{(0, i) \in j} \PP \big( Z^{0} \text{ hits } Z^{i} \big).
\end{equation*}
We fix $j$ and consider $i>0$ so that $(0, i) \in j$ (the case $i < 0$ is similar). That is the
interval of times $s$ with
\begin{equation*}
(s-i)/i  \in (9/10,  11/10), \,\,\, s \in [i+2^{n_{j}} , i+2^{n_{j}+1}]
\end{equation*}
is at least $2$ in length: we note that for each
$i \in \left( \frac{5}{4} 2^{n_j}, \frac{7}{4} 2^{n_j} \right )$ the relevant interval,
$\frac{19}{10} i \leq s \leq \frac{21}{10} i$ is an interval of length greater than
$\frac{5}{4} 2^{n_{j}} \frac{1}{5} = \frac{2^{n_{j}}}{4}$.

\begin{lemma}
\label{lem11}
There exists $c_{2} \in (0, \infty)$ so that for any interval $I$ of length at least $1$
contained in $(1, \infty)$,
\begin{equation*}
\frac{1}{c_{2}} \int_{I} p_t(0,0) \,dt
\leq \PP\big(X^0(t) =0 \text{ for some } t \in I\big)
\leq c_{2} \int_{I} p_t(0,0) \,dt
\end{equation*}
\end{lemma}
\begin{proof}
Consider random variable $W =\int_{a}^{b+1} 1_{\{X^0(s) =0\}}\, ds$ for $I= [a,b]$.
Then
\begin{equation*}
\EE(W) = \int_{a}^{b+1} \PP(X^0(s) =0)\, ds
= \int_{a}^{b+1} p_{s}(0,0)\, ds
\leq 2 \int_{I} p_{s}(0,0) \,ds \, ,
\end{equation*}
by monotonicity of $p_{s}(0,0)$ and the fact that $b-a \geq 1$. But for
$\tau := \inf \{s \in I\colon X^0(s) =0 \}$ we have
$\EE(W \vert {{\cal{F}}}_{ \tau}) \geq e^{-1}$ on $\{\tau < \infty\}$ so
\begin{equation*}
\PP(\tau < \infty) = \PP\big(X^0(t) = 0 \text{ for some } t \in I\big) \leq \EE (W) e
\leq 2 e \int_{I} p_{s}(0,0) ds \, .
\end{equation*}
Equally for $W^\prime = \int_{a}^{b} 1_{\{X^0(s) = 0\}} \, ds$, we have
\begin{equation*}
\EE(W^\prime \,|\, {{\cal{F}}}_{\tau}) \leq \gamma = \int_{0}^{\infty} p_{s} (0,0) \, ds
\quad \text{on } \{\tau < \infty\}
\end{equation*}
and so
\begin{equation*}
\PP(\tau < \infty) \geq \frac{\EE(W^\prime)}{\gamma} =  \frac{1}{\gamma} 
\int_{a}^{b} p_{s}(0,0) \, ds \, .
\end{equation*}
\qed
\end{proof}
\begin{proposition}
\label{}
For some universal $c_{3} \in (0, \infty)$,
\begin{equation*}
c_{3}^{-1} 2^{2n_{j}} p_{2^{n_{j}}}(0,0)
\leq \sum_{(0, i) \in j} \PP\big(Z^{0}\text{ hits }Z^{i}\big)
\leq c_{3} 2^{2n_{j}} p_{2^{n_{j}}} (0,0)\, .
\end{equation*}
\end{proposition}
\begin{proof}
We consider first the upper bound. There are less than $2^{n_{j}}$ relevant $i$.
For such an $i$,
\begin{equation*}
\PP\big(Z^{0}\text{ hits }Z^{i}\big) \leq  \PP \big(X^0(t) \text{ hits } 0 \text{ for some } 
t\in \big[a+2^{n_{j}}, a+3 \cdot 2^{n_{j}}\big]\big)
\end{equation*}
for some $a\geq 0$. By monotonicity of $t \to p_t(0,0)$ and using Lemma \ref{lem11}, 
this is bounded by
\begin{equation*}
c_{2} \int_{2^{n_{j}}}^{3 \cdot 2^{n_{j}}} p_{s}(0,0) \, ds
\leq c_{3} 2^{n_{j}}p_{2^{n_{j}}}(0,0)
\end{equation*}
for some $c_3>0$. On the other side the number of 
$i \in \left( \frac{5}{4} 2^{n_j} , \frac{7}{4} 2^{n_j} \right)$
is greater than $c_{1} \frac{1}{3} 2^{n_j}$ if $R$ was fixed
sufficiently large and for each such $i$, $(\frac{9}{10}i$,
$\frac{11}{10} i) \subset [2^{n_{j}}, 2^{n_{j}+1}]$. Moreover, we have
\begin{eqnarray*}
\PP\big(Z^{0}\text{ hits }Z^{i}\big) &\geq& \frac{1}{c_2} \, \int_{\frac{28}{10}i}^{\frac{32}{10} i} p_{s}(0,0) \, ds
\geq \frac{1}{c_2} \, \frac{4}{10}\,i p_{2^{n_j+3}}(0,0)
\geq \frac{1}{c_2} \, 2^{n_{j}-1} p_{2^{n_j+3}} (0,0)\\
&\geq& \frac{2^{-13}}{c_2} \, 2^{n_{j}}p_{2^{n_{j}-1}}(0,0)
\geq c_3^{-1} \, 2^{n_{j}}p_{2^{n_{j}}}(0,0),
\end{eqnarray*}
because of Lemma \ref{lem11}, Corollary \ref{lem5.1} (by our choice of $j$), monotonicity of $t \to p_t(0,0)$ and possibly after increasing $c_3$.
\qed
\end{proof}
Thus, using that $j\in[k_{0},k_{1}]$ (recall \ref{k1def}), we have a universal $c_{4}$ such that
\begin{equation*}
c_{4} \alpha \geq \sum_{j} \sum_{(0,i)\in j}
\PP\big(Z^{0}\text{ hits }Z^{i} \big)
\geq \frac{\alpha}{c_{4}}.
\end{equation*}
There are two issues to address
\begin{description}[(a)]
\item[(a)]
to show that
\begin{equation*}
\PP\big(\exists j, \exists i \text{ so that } (0,i) \in j, \,\,\, Z^{0}\text{ hits }Z^{i}\big)
\end{equation*} 
is of the order $\alpha$;
\item[(b)]
to show that (a) holds with $Z^{0}$, $Z^{i}$ replaced by our coalescing random
walks $Y^{0}$, $Y^{i}$.
\end{description}
In fact both parts are resolved by the same calculation.

We consider the probability that random walk $Z^{0}$ is involved in a ``3-way''
collision with $Z^{i}$ and $Z^{i^\prime}$ either due to $Z^{0}$ hitting $Z^{i}$ in the
appropriate time interval and then hitting $Z^{i^\prime}$, or $Z^{0}$ hitting $Z^{i}$ 
and, subsequently $Z^{i}$ hitting $Z^{i^\prime}$.
The first case is important to bound so that one can use simple Bonferroni bounds
to get a lower bound on $\PP( \exists i \text{ so that } Z^{i} \text{ hits }Z^{0} )$.
The second is to take account of the fact that we are interested in the future
coalescence of a given random walk $Y^{0}$. As already noted, we can couple
the systems in the usual way so that for all $t$,
$\cup_{i} \{Y_{t}^{i} \} \subseteq \cup_{i} \{Z_{t}^{i}\}$.
The problem is that if for some $i$,  $Z^{i}$ hits $Z^{0}$ due to coalescence this
need not imply that $Y^{i}$ hits $Y^{0}$: if the $Y^{i}$ particles coalesced with a
$Y^{i^\prime}$ before $Z^{i}$ hits $Z^{0}$. Fortunately this event is contained in the union
of events above over $i$, $i^\prime$.

\begin{proposition}
\label{bte}
There exists universal constant $K$ so that for all  $i$ and $i^\prime$ with 
$(0,i^\prime) \in j^\prime$
\begin{equation*}
\PP\big(Z^{0} \text{ hits }Z^{i} \text{ and then }Z^{i^\prime}\big)
\leq K  2^{n_{j^\prime}}p_{2^{n_{j^\prime}}}(0,0) 
\PP\big(Z^{0} \text{ hits }Z^{i}\big)\,.
\end{equation*}
\end{proposition}

\begin{proof}
There are several cases to consider:
$i<0<i^\prime$, $i<i^\prime<0$, $i^\prime<i<0$,  $i^\prime<0<i$, $0<i<i^\prime$ and 
$0<i^\prime<i$.
All are essentially the same so we consider explicitly $0<i<i^\prime$. We leave the reader
to verify that the other cases are analogous.
We choose $j$, $j^\prime$ so that $(0, i) \in j$ and $(0, i^\prime) \in j^\prime$ 
(so necessarilly $j^\prime \geq j$).
We condition on $T_{j}$, $Z^{i}(T_{j})(= Z^{0}(T_{j}))$, for
\begin{equation*}
T_{j} := \inf \bigg\{s \in \left(\frac{19i}{10}, \frac{21i}{10} \right) \cap 
\Big[i+2^{n_{j}}, i+2^{n_{j}+1}\Big] \colon Z^{i}(s) = Z^{0}(s)\bigg\} < \infty.
\end{equation*}
With $x=Z^0(T_{j})$ we have
\begin{eqnarray*}
&&\PP\left( \exists s^\prime \geq T_j \in \left(\frac{19i^\prime}{10}, \frac{21i^\prime}{10} \right) \cap 
\Big[i^\prime+2^{n_{j^\prime}}, i^\prime+2^{n_{j^\prime}+1}\Big]
\colon Z^{i^\prime} (s^\prime)= Z^{0}(s^\prime) \,\,\,\Big\vert\,\,\, G^{0,i}\right)\\
&&\qquad=
 \PP \big(Z^0(t) = x  \mbox{ for some } t \in I_j\big),
\end{eqnarray*}
where $I_j$ is the image of the interval 
\begin{equation*}
\Big[(i^\prime+2^{n_{j^\prime}}) \vee T_{j} \vee \frac{19i^\prime}{10}, 
i^\prime+2^{n_{j^\prime}+1} \wedge \frac{21i^\prime}{10}\Big] \bigg),
\end{equation*}
by the function $t \mapsto 2t-T_j-i^\prime$, for $G^{0, i} = \sigma (Z^{0}(s), Z^{i}(s) \colon s \leq T_{j})$.
By elementary algebra this is less than
\begin{equation*}
\PP \left( Z^0(t) = x \text{ for }t \in \left(\frac{9i^\prime}{10}, \frac{16i^\prime}{5} \right)\right),
\end{equation*}
but by arguing as in  Lemma $\ref{lem11}$, this is bounded by
\begin{eqnarray*}
c_2 \int^{\frac{16i^\prime}{5}}_{\frac{9i^\prime}{10}}p_{s}(0,x) \,ds
&\leq& c_2\int^{\frac{16}{5}i^\prime}_{\frac{9}{10}i^\prime} p_{s}(0,0) \,ds 
\leq c_2 \frac{23}{10}\,i^\prime p_{2^{n_{j^\prime}-1}}(0,0) \\
&\leq& c_2 \frac{23}{9} 2^{13}2^{n_{j^\prime}}p_{2^{n_{j^\prime}+3}}(0,0) 
\leq c^{\prime} 2^{n_{j^\prime}}p_{2^{n_{j^\prime}}}(0,0)
\end{eqnarray*}
for some universal constant $c^\prime$, where we use 
symmetry and monotonicity of $p_{s}(\cdot,\cdot)$ and Corollary \ref{lem5.1}, by the choice of our $n_{j^\prime}$.
So given that
\begin{equation*}
\PP(T_{i} < \infty) = \PP\big(Z^{0} \text{ hits }Z^{i} \big),
\end{equation*}
the desired bound is achieved.
\qed
\end{proof}

\begin{corollary}
\label{cor12}
For $\alpha$ sufficiently small
\begin{equation*}
\PP\big( \exists i \colon Z^{0}\text{ hits }Z^{i}\big) \geq \frac{\alpha}{2}\, .
\end{equation*}
\end{corollary}

\begin{proof}
By Bonferroni, the desired probability is superior to
\begin{equation*}
\sum_{i} \PP\big(Z^{0}\text{ hits }Z^{i}\big)
- \sum_{i, i^\prime} \PP\big(Z^{0}\text{ hits }Z^{i}\text{ and then }Z^{i^\prime}\big)
\geq \alpha - Kc_{3}^2\alpha^{2} \geq \frac{\alpha}{2}
\end{equation*}
if $\alpha \leq 1/(2Kc_3^2)$.
\qed
\end{proof}

We similarly show
\begin{proposition}
\label{prop7}
\noindent There exists universal constant $K$ so that for all $i^\prime$ and 
$i$ with $(0,i)  \in j$,
\begin{equation*}
\PP\big(Z^{0}\text{ hits }Z^{i}\text{ after }Z^{i}\text{ hits }Z^{i^\prime}\big)
\leq K 2^{n_{j}}p_{2^{n_{j}}}(0,0) \PP\big(Z^{i}\text{ hits }Z^{i^\prime}\big)\, .
\end{equation*}
\end{proposition}
This gives as a corollary
\begin{corollary}
\label{cor13}
For the coalescing system $\{Y^{i}\}_{i\in V}$ provided $\alpha$ is sufficiently small,
\begin{equation*}
\PP\big(Y^{i}\text{ dies after time }R\big) \geq \frac{\alpha}{5}\, .
\end{equation*}
\end{corollary}

\begin{proof}
We have of course from the labeling scheme
\begin{eqnarray*}
&&\PP\big(Y^{i} \text{ dies after time }R\big)\\
&&\qquad\geq\ \frac{1}{2} \PP\big(Y^{i} \text{ hits }Y^{i^\prime}\text{ in appropriate 
time interval for some }i^\prime\big)\\
&&\qquad\geq \frac{1}{2} \PP\big(Z^{i}\text{ hits }Z^{i^\prime}\text{ in appropriate 
time interval for some }i^\prime\big) \\
&&\qquad- \frac{1}{2} \PP\big(Z^{i}\text{ hits }Z^{i^\prime}\text{ in appropriate 
time interval for some }i^\prime \text{ so that }\\
&&\qquad\qquad\qquad Z^{i^\prime}\text{ hits some }Z^{i^{\prime\prime}}
\text{ previously}\big)\\
&&\qquad\geq \frac{\alpha}{4} - Kc_3^2 \alpha^{2} \geq \frac{\alpha}{5}
\end{eqnarray*}
for $\alpha \leq 1/(20Kc_3 ^2)$.
\qed
\end{proof}
We can now complete the proof of Proposition \ref{prop3} and hence that of
Proposition \ref{prop2}. If we have $c_0 > 0$, then we can find $0<\epsilon< \alpha / 200$
and $\alpha$ so small that the relevant results above hold, in particular Corollary \ref{cor13}.
Thus the density of Y's is reduced by at least $\alpha / 5$.  But by our choice of $\epsilon$ 
and Lemma \ref{lem8}, the density of $X$'s is reduced by at least 
$\alpha/5 -6 \epsilon \geq \alpha / 6 \geq 3 \epsilon $ which is a contradiction with Proposition
\ref{lem3}, because it would entail the density falling strictly below $c_{0}$.


\section{Proof of Theorem \ref{th3}}
\label{S4}

In what follows we assume, as in Section \ref{S2}, that $p=1$, the extension
to arbitrary $p\geq 1$ being straightforward.

We begin by specifying the random walk $(X(t))_{t\geq 0}$ on $\Z^4$ defined by
\begin{equation*}
\label{rw1.4}
X(t) = S(t)+ e_{1}N(t)
\end{equation*}
with $(S(t))_{t\geq 0}$ denoting a simple random walk on $\Z^4$,
$(N(t))_{t\geq 0}$ a rate $1$ Poisson process and $e_{1}=(1,0,0,0)$ the first unit vector
in $\Z^4$.

Thus our random walk $(X(t))_{t\geq 0}$ is highly transient but its symmetrization
is a mean zero random walk and by the local central limit theorem, we have
\begin{equation*}
\label{}
\int_{0}^{\infty} t p_{t}^{(s)}(0,0) \,dt = \infty,
\end{equation*}
where $p_{t}(\cdot,\cdot)$ is the semigroup associated to $(X(t))_{t\geq 0}$.

It remains to show that $\lambda_{p}(\kappa)<\gamma$ for all $\kappa\in[0,\infty)$.
Our approach is modeled on the proof of the first part of Theorem \ref{th1}. We wish
again to pick \emph{bad environment set} $B_{E}$ associated to the $\xi$-process and
\emph{bad random walk set} $B_{W}$ associated to the random walk $X^{\kappa}$
so that
\begin{eqnarray}
\label{lambdaupbd3}
&&\EE
\Big(\exp\Big[\gamma\int_{0}^n\xi_{s}(X^\kappa(s))\, ds\Big]\Big)\\
&&\qquad\leq \Big(\PP(B_{E})+\PS(B_{W})\Big)e^{\gamma n}
+ \EE\Big(\one_{B_{E}^\c\cap B_{W}^\c}
\exp\Big[\gamma\int_{0}^n\xi_{s}(X^\kappa(s))\, ds\Big]\Big)
\nonumber\\
\end{eqnarray}
with, for some $0<\delta<1$,
\begin{equation}
\label{badbound}
\PP(B_{E})\leq e^{-\delta n},
\qquad
\PS(B_{W})\leq e^{-\delta n},
\end{equation}
and, automatically from the definition of $B_{E}$ and $B_{W}$,
\begin{equation}
\label{goodbound}
\int_{0}^{n} \xi_{s}(X^\kappa(s))\,ds \leq n(1-\delta)
\qquad\text{on } B_{E}^\c\cap B_{W}^\c
\end{equation}
(as in the proof of Theorem \ref{th1}).
Since, combining (\ref{lambdaupbd3}--\ref{goodbound}), we obtain
\begin{equation*}
\label{}
\lim_{n\to\infty}\frac{1}{n}\log
\EE\Big(\exp\Big[\gamma\int_{0}^n\xi_{s}(X^\kappa(s))\, ds\Big]\Big)
< \gamma\, ,
\end{equation*}
it is enough to prove (\ref{badbound}). All of this has been done in the proof
of Theorem \ref{th1} in a different situation. The major difference is that we
need to modify the collection of skeletons used.
\begin{lemma}
\label{lem1-countex} Let $X^{\kappa} (\cdot)$ be a speed $\kappa$ simple
random walk in four dimensions. Fix $M\in\N\setminus \{1\}$. There exists
$c>0$ so that for $M$ large and all $n$,
outside of an $e^{-cn}$ probability event, there exists
$0\leq i_{1}<i_{2}<\cdots < i_{n/2M}\leq n$ so that
\begin{equation}
\label{driftSkl}
X^{\kappa}_{(1)}(i_{j}M+kM) - X^{\kappa}_{(1)}(i_{j}M) > - \frac{kM}{2} ,
\qquad  j\in\{1,\cdots,n/(2M)\},\, k\geq 0\, ,
\end{equation}
where $(X^{\kappa}_{(1)}(t))_{t\geq 0}$ denotes the first coordinate of
$(X^{\kappa}(t))_{t\geq 0}$.
\end{lemma}
\begin{proof}
Define
\begin{equation*}
\sigma_{1} = \inf \big\{kM >0\colon X^{\kappa}_{(1)}(kM) \leq -kM/2 \big\}
\end{equation*}
and recursively
\begin{equation*}
\sigma_{i+1} = \inf \big\{kM >\sigma_{i}\colon
X^{\kappa}_{(1)}(kM)-X_{(1)}^\kappa(\sigma_{i})
\leq -(kM-\sigma_{i})/2 \big\} \, .
\end{equation*}
Since the event
\begin{equation*}
\{rM \leq \sigma_{1} < \infty \}  \subset  \
\cup_{k=r}^{\infty }\big\{X^{\kappa}_{(1)}(kM) \leq  -kM/2 \big\} \, ,
\end{equation*}
we have  easily that, for all $r$,
\begin{equation*}
\PP (rM \leq \sigma_{1} < \infty) \ \leq \  e^{-rMc}
\end{equation*}
for $c>0$ not depending on $n$ or $M$.
If we now define
\begin{equation*}
\tau_{1} = \inf \big\{kM >0\colon X_{(1)}^\kappa(jM)-X^{\kappa}_{(1)}(kM) > -(j-k)M/2
\,\,\,\forall j>k\big\}
\end{equation*}
and recursively
\begin{equation*}
\tau_{i+1} = \inf \big\{kM >\tau_{i}\colon
X^{\kappa}_{(1)}(jM)-X^{\kappa}_{(1)}(kM)  > -(j-k)M/2
\,\,\,\forall j>k \big\}\, ,
\end{equation*}
it is easily seen that
\begin{eqnarray*}
\PP(\tau_{1}  \geq rM)
&\leq& \PP(\exists 1\leq k \leq r\colon rM \leq \sigma_{k} < \infty)\\
&\leq&  \sum_{k=1}^{r-1}\,\,\,  \sum_{0<x_{1 }< \cdots < x_{k} < r}
\PP\big(\sigma_{i} =x_{i}M \,\,\,\forall i \leq k, \, rM \leq \sigma_{k+1} < \infty \big)\\
&\leq& \sum_{k=1}^{r-1}\,\,\,  \sum_{0<x_{1 }< \cdots < x_{k} <r} e^{-rMc} 
\leq \ e^{-rMc}2^r
\end{eqnarray*}
which is less than $e^{-rMc/2}$ if $M$ is fixed sufficiently large.
We have
\begin{itemize}
\item
$(\tau_{i+1} - \tau_{i})_{i \geq 1}$ are i.i.d. (this follows from Kuczek's argument
(see \cite{kuc89})).
\item
Provided $M$ has been fixed sufficiently large for each integer
$r \geq 1$, $\PP(\tau_{i+1} - \tau_{i} \geq rM) \leq e^{-rMc/4}$.
This follows from the fact that random variable $\tau_{i+1} - \tau_{i}$
is simply the random variable $\tau_{1}$ conditioned on an event of
probability at least $1/2$ (provided $M$ was fixed large).
\end{itemize}
Thus by elementary properties of geometric random variables we have
\begin{equation*}
\PP(\tau_{n/2M} >n) \leq \PP\left(Y \geq \frac{n}{2M}\right) \leq e^{-cn}
\end{equation*}
for $Y \sim B\left(\frac{n}{M},e^{-cM/4}\right)$, $c>0$ and $M$ large. 
This completes the proof of the lemma.
\qed
\end{proof}

Given that the path of the random walk satisfies the condition of this lemma,
we call the (not uniquely defined) points $i_{1}, i_{2}, \cdots $ \emph{regular } points.

Given this result, we consider the $M$-skeleton induced by the values
$X^\kappa (jM)$, $0\leq j \leq n/M\}$, discretized via spatial cubes of length $M/8$
(rather than $2M$ as in the proof of Theorem \ref{th1}). It is to be noted
that if $(X^ \kappa (t))_{0\leq t\leq n}$ satisfies the claim for Lemma \ref{lem1-countex}
and $y_{0}:=0, y_{1}, y_{2}, \cdots, y_{n/M}$ with $y_{k}\in \Z^4$,
$0\leq k\leq n/M$, is its $M$-skeleton, namely,
\begin{equation}
\label{Xkappaskl*}
X^\kappa(kM)\in C_{y_{k}}
:=\prod_{j=1}^{4}\bigg[y_{k}^{(j)} \frac{M}{8}, (y_{k}^{(j)}+1)\frac{M}{8}\bigg),
\qquad 0\leq k\leq n/M,
\end{equation}
where $y_{k}^{(j)}$ denotes the $j$-th coordinate of $y_{k}$ (we suppose
without loss of generality that $M$ is a multiple of $8$). Then, by (\ref{driftSkl})
and (\ref{Xkappaskl*}), we must have
\begin{equation*}
y_{i_{j}}^{(1)}-4k\leq y_{i_{j}+k}^{(1)} +1.
\end{equation*}
In particular, we must have
\begin{equation}
\label{sklcond}
y^{(1)}_{i_{j^\prime}} -4 (i_{j}-i_{j^\prime})
\leq y^{(1)}_{i_{j}} +1
\qquad \forall i_{j^\prime}<i_{j}.
\end{equation}

In the following we modify the definition of appropriate skeletons by adding in the
requirement that the skeleton must possess at least $n/2M$ indices
$i_{1}, i_{2}, \ldots, i_{n/2M}$ with the corresponding $y_{i_{j}}$ satisfying (\ref{sklcond}).
We note that the resizing of the cubes makes the notion of acceptability a little more
stringent but does not change the essentials.

Remark first that Lemma \ref{bw2} is still valid in our new setting. Lemma \ref{lem1-countex} 
immediately gives that with this new definition, Lemma \ref{bw1} remains true.
Of course since this definition is more restrictive we have
\begin{equation*}
\label{}
|\Xi_{A}|\leq K^{n/M}
\end{equation*}
for $K$ as in Lemma \ref{appropriate}. 

In fact in our program all that remains to do, that is in any substantive way
different from the proof of Theorem \ref{th1}, is to give a bound on the probability
of $B_{E}$ for appropriate $B_{E}$.  This is the content of
the lemma below (analogous to Lemma \ref{be.lem}).
Given this lemma, we can then proceed exactly as with the proof of Theorem \ref{th1}.
\begin{lemma}
For any skeleton $(y_{k})_{0\leq k\leq n/M}$ in $\Xi_{A}$, the probability that $\xi $ is
not good for $(y_{k})_{0\leq k\leq n/M}$,
i.e.,
\begin{equation*}
\not\exists \frac{n}{4M} \text{ indices } 1\leq j\leq \frac{n}{M} \colon \xi_{s}(z)=0\,\,\,
\forall s\in[jM,jM+1]
\text{ for some } z\in C_{y_{j}}\, ,
\end{equation*}
is less than $(4K)^{-n/M}$.
\end{lemma}
\begin{proof}
We note that proving the analogous result for Theorem \ref{th1}, we did not need
our skeleton to be in $\Xi_{A}$, the proof worked over any skeleton. For us however it is vital
that our skeleton satisfies (\ref{sklcond}).

We consider a skeleton in $\Xi_{A}$. Let the first $n/2M$ regular points of our skeleton
be $ i_{1}, i_{2}, \cdots i_{n/2M}$.  
For each $1\leq i_{j}\leq n/2M$, we choose $R$ points
\begin{equation*}
x_{1}^{i_{j}},\ldots, x_{R}^{i_{j}} \in C_{ y_{i_{j}}}
\end{equation*}
so spread out that  for random walks $(X(t))_{t\geq 0}$ as in (\ref{rw1.4}) beginning at the points 
$x_{k}^{i_{j}}$, $k=1,\cdots,R$, the chance that two of them meet is less than 
$0<\epsilon\ll 1$.

Now, consider $i_{j^\prime}<i_{j}$ and the probability that a random walk starting at
$(x_{k}^{i_{j}}, i_{j} M)$ meets a random walk starting at 
$(x_{k^\prime}^{i_{j^\prime}}, i_{j^\prime} M)$ satisfies the following lemma.
\begin{lemma}
For $i_{j^\prime}<i_{j}$, there exits $K>0$ such that
\begin{equation*}
\PP\bigg(X^{x_{k}^{i_{j}},i_{j} M} \,\,\text{ meets }\,\, 
X^{x_{k^\prime}^{i_{j^\prime}},i_{j^\prime} M}\bigg)\leq \frac{K}{M^2(i_{j}-i_{j^\prime})^2} \, .
\end{equation*}
\end{lemma}
\begin{proof}
The important point is that since our skeleton is in $\Xi_{A}$,
\begin{equation}
\label{imp.lem4.3}
\Big(x_{k^\prime}^{i_{j^\prime}}\Big)^{(1)} \leq \Big(x_{k}^{i_{j}}\Big)^{(1)} 
+ (i_{j}-i_{j^\prime})\frac{M}{2}+ \frac{M}{4}\, ,
\end{equation}
and so we have
\begin{eqnarray*}
&&\PP\bigg(X^{x_{k}^{i_{j}},i_{j} M} \,\,\text{ meets }\,\,
X^{x_{k^\prime}^{i_{j^\prime}},i_{j^\prime} M}\bigg) \\
&&\qquad\leq \PP\bigg(\Big(X^{x_{k}^{i_{j}},i_{j} M}\Big)^{(1)}
\left((i_{j}-i_{j^\prime})M\right)
\geq \big(x_{k}^{i_{j}}\big)^{(1)} + (i_{j}-i_{j^\prime})\frac{3M}{4}\, ,\\
&& \qquad\qquad\quad X^{x_{k}^{i_{j}},i_{j} M} \,\,\text{ meets }\,\,
X^{x_{k^\prime}^{i_{j^\prime}},i_{j^\prime} M}\bigg) \\
&&\qquad + \PP\bigg(\Big(X^{x_{k}^{i_{j}},i_{j} M}\Big)^{(1)}
\left((i_{j}-i_{j^\prime})M\right)
\leq \big(x_{k}^{i_{j}}\big)^{(1)} + (i_{j}-i_{j^\prime})\frac{3M}{4}\bigg).
\end{eqnarray*}
By standard large deviations bounds,
\begin{equation*}
\PP\bigg(\Big(X^{x_{k}^{i_{j}},i_{j} M}\Big)^{(1)}
\left((i_{j}-i_{j^\prime})M\right)
\leq \big(x_{k}^{i_{j}}\big)^{(1)} + (i_{j}-i_{j^\prime})\frac{3M}{4}\bigg)
\leq e^{-CM(i_{j}-i_{j^\prime})},
\end{equation*}
for some universal $C\in (0,\infty)$. For the other term, we have the following lemma.
\begin{lemma}
For two independent processes $X=(X(t))_{t\geq 0}$ and $Y=(Y(t))_{t\geq 0}$
with $X(0)=x\in\Z^4$ and $Y(0)=y\in\Z^4$, the probability that $X$ ever meets
$Y$ is bounded by $K/\|x-y\|_{\infty}^2$.
\end{lemma}
\begin{proof}
$X-Y$ is not exactly a simple random walk, but it is a symmetric random walk
and so local CLT gives appropriate random walks bounds (see, e.g., \cite{law10}).
\qed
\end{proof}
From this and inequality (\ref{imp.lem4.3}), we have
\begin{eqnarray*}
&&\PP\bigg(X^{x_{k}^{i_{j}},i_{j} M} \,\text{meets }
X^{x_{k^\prime}^{i_{j^\prime}},i_{j^\prime} M}  \Big\vert
\Big(X^{x_{k}^{i_{j}},i_{j} M}\Big)^{(1)}
\left((i_{j}-i_{j^\prime})M\right)
\geq \big(x_{k}^{i_{j}}\big)^{(1)} + (i_{j}-i_{j^\prime})\frac{3M}{4}\bigg)\\
&&\qquad\leq \frac{K}{M^2 (i_{j}-i_{j^\prime})^2}\, .
\end{eqnarray*}
\qed
\end{proof}
Thus, for any $R$ large but fixed, we can choose $M$ so that for all skeletons
in $\Xi_{A}$ and each $i_{j^\prime}$, we have
\begin{eqnarray*}
\sum_{i_{j}<i_{j^\prime}}\sum_{k,k^\prime}
\PP\bigg(X^{x_{k}^{i_{j}},i_{j} M} \,\,\text{ meets }\,\,
X^{x_{k^\prime}^{i_{j^\prime}},i_{j^\prime} M}\bigg)
&\leq& \frac{R^2 K}{M^2}\sum_{r=1}^{+\infty}\frac{1}{r^2}\\
&\leq& \frac{R^2 K^\prime}{M^2} < \epsilon^2
\end{eqnarray*}
with $M$ chosen sufficiently large, which is analogous to (\ref{2.2.08}).
From this point on, the rest follows as for the proof of Lemma \ref{be.lem}.
\qed
\end{proof}
\input{referenc}

\end{document}

%% file: referenc.tex
%
%
%